\magnification\magstep1

\font\ninerm=cmr9  \font\eightrm=cmr8  \font\sixrm=cmr6
\font\ninei=cmmi9  \font\eighti=cmmi8  \font\sixi=cmmi6
\font\ninesy=cmsy9 \font\eightsy=cmsy8 \font\sixsy=cmsy6
\font\ninebf=cmbx9 \font\eightbf=cmbx8 \font\sixbf=cmbx6
\font\nineit=cmti9 \font\eightit=cmti8 
\font\ninett=cmtt9 \font\eighttt=cmtt8 
\font\ninesl=cmsl9 \font\eightsl=cmsl8

\font\twelverm=cmr12 at 15pt
\font\twelvei=cmmi12 at 15pt
\font\twelvesy=cmsy10 at 15pt
\font\twelvebf=cmbx12 at 15pt
\font\twelveit=cmti12 at 15pt
\font\twelvett=cmtt12 at 15pt
\font\twelvesl=cmsl12 at 15pt
\font\twelvegoth=eufm10 at 15pt

\font\tengoth=eufm10  \font\ninegoth=eufm9
\font\eightgoth=eufm8 \font\sevengoth=eufm7 
\font\sixgoth=eufm6   \font\fivegoth=eufm5
\newfam\gothfam \def\goth{\fam\gothfam\tengoth} 
\textfont\gothfam=\tengoth
\scriptfont\gothfam=\sevengoth 
\scriptscriptfont\gothfam=\fivegoth

\catcode`@=11
\newskip\ttglue

\def\tenpoint{\def\rm{\fam0\tenrm}
  \textfont0=\tenrm \scriptfont0=\sevenrm
  \scriptscriptfont0\fiverm
  \textfont1=\teni \scriptfont1=\seveni
  \scriptscriptfont1\fivei 
  \textfont2=\tensy \scriptfont2=\sevensy
  \scriptscriptfont2\fivesy 
  \textfont3=\tenex \scriptfont3=\tenex
  \scriptscriptfont3\tenex 
  \textfont\itfam=\tenit\def\it{\fam\itfam\tenit}%
  \textfont\slfam=\tensl\def\sl{\fam\slfam\tensl}%
  
  \textfont\ttfam=\tentt\def\tt{\fam\ttfam\tentt}%
  \textfont\gothfam=\tengoth\scriptfont\gothfam=\sevengoth 
  \scriptscriptfont\gothfam=\fivegoth
  \def\goth{\fam\gothfam\tengoth}
  \textfont\bffam=\tenbf\scriptfont\bffam=\sevenbf
  \scriptscriptfont\bffam=\fivebf
  \def\bf{\fam\bffam\tenbf}%
  \tt\ttglue=.5em plus.25em minus.15em
  \normalbaselineskip=12pt \setbox\strutbox\hbox{\vrule
  height8.5pt depth3.5pt width0pt}%
  \let\big=\tenbig\normalbaselines\rm}

\def\ninepoint{\def\rm{\fam0\ninerm}
  \textfont0=\ninerm \scriptfont0=\sixrm
  \scriptscriptfont0\fiverm
  \textfont1=\ninei \scriptfont1=\sixi
  \scriptscriptfont1\fivei 
  \textfont2=\ninesy \scriptfont2=\sixsy
  \scriptscriptfont2\fivesy 
  \textfont3=\tenex \scriptfont3=\tenex
  \scriptscriptfont3\tenex 
  \textfont\itfam=\nineit\def\it{\fam\itfam\nineit}%
  \textfont\slfam=\ninesl\def\sl{\fam\slfam\ninesl}%
  \textfont\ttfam=\ninett\def\tt{\fam\ttfam\ninett}%
  \textfont\gothfam=\ninegoth\scriptfont\gothfam=\sixgoth 
  \scriptscriptfont\gothfam=\fivegoth
  \def\goth{\fam\gothfam\tengoth}
  \textfont\bffam=\ninebf\scriptfont\bffam=\sixbf
  \scriptscriptfont\bffam=\fivebf
  \def\bf{\fam\bffam\ninebf}%
  \tt\ttglue=.5em plus.25em minus.15em
  \normalbaselineskip=11pt \setbox\strutbox\hbox{\vrule
  height8pt depth3pt width0pt}%
  \let\big=\ninebig\normalbaselines\rm}

\def\eightpoint{\def\rm{\fam0\eightrm}
  \textfont0=\eightrm \scriptfont0=\sixrm
  \scriptscriptfont0\fiverm
  \textfont1=\eighti \scriptfont1=\sixi
  \scriptscriptfont1\fivei 
  \textfont2=\eightsy \scriptfont2=\sixsy
  \scriptscriptfont2\fivesy 
  \textfont3=\tenex \scriptfont3=\tenex
  \scriptscriptfont3\tenex 
  \textfont\itfam=\eightit\def\it{\fam\itfam\eightit}%
  \textfont\slfam=\eightsl\def\sl{\fam\slfam\eightsl}%
  \textfont\ttfam=\eighttt\def\tt{\fam\ttfam\eighttt}%
  \textfont\gothfam=\eightgoth\scriptfont\gothfam=\sixgoth 
  \scriptscriptfont\gothfam=\fivegoth
  \def\goth{\fam\gothfam\tengoth}
  \textfont\bffam=\eightbf\scriptfont\bffam=\sixbf
  \scriptscriptfont\bffam=\fivebf
  \def\bf{\fam\bffam\eightbf}%
  \tt\ttglue=.5em plus.25em minus.15em
  \normalbaselineskip=9pt \setbox\strutbox\hbox{\vrule
  height7pt depth2pt width0pt}%
  \let\big=\eightbig\normalbaselines\rm}

\def\twelvepoint{\def\rm{\fam0\twelverm}
  \textfont0=\twelverm\scriptfont0=\tenrm
  \scriptscriptfont0\sevenrm
  \textfont1=\twelvei\scriptfont1=\teni
  \scriptscriptfont1\seveni 
  \textfont2=\twelvesy\scriptfont2=\tensy
  \scriptscriptfont2\sevensy 
   \textfont\itfam=\twelveit\def\it{\fam\itfam\twelveit}%
  \textfont\slfam=\twelvesl\def\sl{\fam\slfam\twelvesl}%
  \textfont\ttfam=\twelvett\def\tt{\fam\ttfam\twelvett}%
  \textfont\gothfam=\twelvegoth\scriptfont\gothfam=\ninegoth 
  \scriptscriptfont\gothfam=\sevengoth
  \def\goth{\fam\gothfam\twelvegoth}
  \textfont\bffam=\twelvebf\scriptfont\bffam=\ninebf
  \scriptscriptfont\bffam=\sevenbf
  \def\bf{\fam\bffam\twelvebf}%
  \tt\ttglue=.5em plus.25em minus.15em
  \normalbaselineskip=12pt \setbox\strutbox\hbox{\vrule
  height9pt depth4pt width0pt}%
  \let\big=\twelvebig\normalbaselines\rm}

\font\bigtitlefont=cmbx12 at 16pt

\def\zmod#1{\,\,({\rm mod}\,\,#1)}
\def\QQ{{\bf Q}}
\def\ZZ{{\bf Z}}
\def\FF{{\bf F}}
\def\HH{{\bf H}}
\def\RR{{\bf R}}
\def\CC{{\bf C}}
\def\PP{{\bf P}}

\def\zmod#1{\,\,({\rm mod}\,\,#1)}

\newcount\newpen\newpen=50
\def\qqlineA#1 #2 #3 {\line{\kern#1truept\vrule height0.5truept depth0.15truept 
width#2truept\hskip#3truept plus0.1truept\vrule height0.5truept depth0.05truept 
width16truept\hskip#3truept plus0.1truept\vrule height0.5truept depth0.15truept 
width#2truept\kern#1truept}}
\def\qqlineB#1 #2 #3 {\line{\kern#1truept\vrule height0.5truept depth0.15truept 
width#2truept\hskip#3truept plus0.1pt\vrule height0.5truept depth0.15truept 
width#2truept\kern#1truept}}
\def\qqlineC#1 #2 {\line{\hskip#1truept plus0.1pt\vrule height0.5truept 
depth0.15truept width#2truept\hskip#1truept plus0.1pt}}
\font\sanserifeight=cmss8 at 8truept
\font\sanserifeightfive=cmss8 scaled\magstep4
\newbox\boxtorre
\setbox\boxtorre=\vtop{\hsize=16truept\offinterlineskip\parindent=0truept
	\line{\vrule height1truept depth0.5truept width1truept\hfil\vrule 
	height1truept depth0.5truept width1truept\hfil\vrule height1truept 
	depth0.5truept width1truept\hfil\vrule height1truept depth0.5truept 
	width1truept\hfil\vrule height1truept depth0.5truept width1truept}
	\line{\vrule height1truept depth1truept width16truept}
	\vskip1.5truept
	\line{\vrule height1truept depth0truept width16truept}
	\vskip1truept
	\line{\vrule height1truept depth0truept width16truept}
	\vskip1truept
	\line{\vrule height1truept depth0truept width16truept}
	\vskip1truept
	\line{\vrule height1truept depth0truept width16truept}
	\vskip1truept
	\line{\vrule height1truept depth0truept width16truept}
	\vskip1truept
	\line{\vrule height1truept depth0truept width16truept}
	\vskip1truept
	\line{\vrule height1truept depth0truept width16truept}
	\vskip1truept
	\line{\vrule height1truept depth0truept width16truept}
	\vskip1truept
	\line{\vrule height1truept depth0truept width16truept}
	\vskip1truept
	\line{\vrule height1truept depth0truept width16truept}
	\vskip1truept}
\newbox\boxu
\setbox\boxu=\hbox{\vrule height1truept depth24truept width4truept\kern6truept
	\copy\boxtorre\kern6truept\vrule height1truept depth24truept 
	width4truept}
\newbox\boxuni
\setbox\boxuni=\vtop{\baselineskip=-1000pt\lineskip=-0.15truept
\lineskiplimit=0pt\parindent=0truept\hsize=36truept
\line{\copy\boxu\hfil}
\qqlineA 0.0 4.1 5.9 
\qqlineA 0.0 4.2 5.8 
\vskip-0.1truept
\qqlineB 0.1 4.2 27.4 
\qqlineB 0.21 4.2 27.2 
\qqlineA 0.3 4.2 5.5 
\qqlineA 0.41 4.2 5.4 
\qqlineB 0.5 4.3 26.4 
\qqlineB 0.6 4.4 26.0 
\qqlineA 0.7 4.5 4.8 
\qqlineA 0.8 4.7 4.5 
\qqlineB 1.0 4.8 24.4 
\qqlineB 1.2 4.8 24.0 
\qqlineB 1.5 4.8 23.4 
\qqlineB 1.7 4.91 22.8 
\qqlineB 1.9 5.1 22.0 
\qqlineB 2.21 5.2 21.2 
\vskip-0.1truept
\qqlineB 2.5 5.3 20.4 
\qqlineB 2.8 5.5 19.4 
\qqlineB 3.1 5.7 18.4 
\qqlineB 3.4 6.0 17.21 
\qqlineB 3.7 6.3 16.0 
\qqlineB 4.2 6.8 14.0 
\qqlineB 4.6 7.4 12.0 
\qqlineB 5.1 7.9 10.0 
\qqlineB 5.6 9.0 6.8  
\vskip-0.1truept
\qqlineC 6.1 23.8 
\qqlineC 6.7 22.6 
\qqlineC 7.4 21.21 
\qqlineC 8.0 20.0 
\qqlineC 8.6 18.8 
\vskip-0.1truept
\qqlineC 9.6 16.8
\qqlineC 10.91 14.2 
\qqlineC 12.0 12.0 
\qqlineC 13.5 9.0 
\vskip 3.1truept
\centerline{\sanserifeight TOR VERGATA}}
\def\today{\ifcase\month\or January \or February\or
March\or April\or May\or June\or July\or August\or
September\or October\or November\or December\fi
\space\number\day, \number\year}

\newbox\boxprelim
\setbox\boxprelim\vtop to\dp\boxuni{\parindent=0truept\hsize=5truecm\eightpoint
Preliminary version\par \today\par\vfill}
\newbox\boxroma
\setbox\boxroma\vtop
to\dp\boxuni{\hsize=6truecm\parindent=0truept \vglue
.0truecm\sanserifeightfive II Universit\`a degli\par
\smallskip \hskip 1truecm Studi di Roma\par\vfill}

\def\title #1\par{
  \hbox{\copy\boxprelim\hskip 1.2in
\copy\boxuni\hskip 1cm\copy\boxroma}  
  \bigskip\bigskip\bigskip 
  \bigskip\bigskip\bigskip\bigskip\noindent
  {\bigtitlefont#1}\par
  \bigskip\bigskip  \bigskip\bigskip\noindent
Pietro Mercuri
\smallskip\noindent
Ren\'e Schoof
  \bigskip\eightpoint\noindent
  \vbox{
  \hbox{Dipartimento di Matematica}
  \hbox{$\hbox{2}^{\hbox{a}}$ Universit\`a di Roma 
  ``Tor Vergata"}
  \hbox{I-00133 Roma ITALY}
  \hbox{Email: \tt mercuri.ptr@gmail.com}
  \hbox{\qquad\quad\ \tt schoof.rene@gmail.com}}
  \bigskip\bigskip
  \tenpoint
  }
\def\abstract #1\par{\eightpoint\vbox{\noindent
  {\bf Abstract.\ }#1\par}\bigskip\tenpoint
  }
\def\bibliography#1\par{\vskip0pt
  plus.3\vsize\penalty-250\vskip0pt
  plus-.3\vsize\bigskip\vskip\parskip
  \message{Bibliography}\leftline{\bf
  Bibliography}\nobreak\smallskip\noindent
  \ninepoint\frenchspacing#1}
\outer\def\beginsection#1\par{\vskip0pt
  plus.3\vsize\penalty\newpen\newpen=-50\vskip0pt
  plus-.3\vsize\bigskip\vskip\parskip
  \message{#1}\leftline{\bf#1}
  \nobreak\smallskip\noindent}

\def\AL{{1}}
\def\BD{{2}}
\def\BC{{3}}
\def\BA{{4}}
\def\BB{{5}}
\def\BPR{{6}}
\def\BK{{7}}
\def\BU{{8}}
\def\KCON{{9}}
\def\DSED{{10}}
\def\DGGS{{11}}
\def\DOSE{{12}}
\def\DMS{{13}}
\def\LAN{{14}}
\def\BAO{{15}}
\def\LI{{16}}
\def\MAG{{17}}
\def\MDB{{18}}
\def\MAI{{19}}
\def\ME{{20}}
\def\MER{{21}}
\def\SAG{{22}}
\def\SD{{23}}
\def\SE{{24}}
\def\ST{{25}}
\def\WAL{{26}}

\title 
Modular forms invariant under non-split Cartan subgroups

\abstract In this paper we describe a method for computing  a basis 
for the space of weight $2$ cusp forms
invariant under  a non-split Cartan subgroup of prime level $p$. As an application we compute,
for certain small values of $p$, explicit equations over $\QQ$
for the canonical embeddings of the associated modular curves.

\beginsection 1. Introduction.

It is well known how to compute  bases 
for the spaces of cusp forms that are invariant under the modular groups
$\Gamma_0(N)$ or~$\Gamma_1(N)$. Indeed, efficient algorithms to
compute $q$-expansions of eigenforms exist~[\MAG, \SAG] and 
extensive tables  are available online~[\BK, \MDB, \ST].
For other congruence subgroups of ${\rm SL}_2(\ZZ)$ the situation is different. 
While for some groups, like {\it split} Cartan subgroups, there are efficient algorithms~[\SAG] and 
it is easy to obtain $q$-expansions from the existing tables for $\Gamma_0(N)$, for other subgroups this is not so immediate~[\BC, \BA].

In this paper we describe a method to compute $q$-expansions of a basis for the space $S_2(\Gamma_{\rm ns}(p))$ of weight $2$ cusp forms invariant under a {\it non-split} Cartan subgroup $\Gamma_{\rm ns}(p)$ of prime level~$p$. 
As in  the computation for $p=13$ by B.~Baran~[\BB], we obtain  a basis of~$S_2(\Gamma_{\rm ns}(p))$
by applying trace maps to certain normalized 
eigenforms in $S_2(\Gamma_0(p^2))$ and~$S_2(\Gamma_1(p))$.
In Baran's computation for $p=13$ this involves,  only one eigenform. It generates a cuspidal ${\rm GL}_2(\FF_p)$-representation.
For larger primes $p$,  several non-isomorphic
irreducible representations, both cuspidal and principal series,  are involved.
This complicates matters, since in each case the trace map is different.
Our main results are the formulas of  Propositions~4.4 and~5.2. 

As an application we are able to compute  
explicit equations  for the canonical embeddings of  the modular curves $X_{\rm ns}(p)$
associated to  the non-split Cartan subgroups and the curves $X^+_{\rm ns}(p)$ associated to their normalizers.
Since our method allows us to compute a basis that is  defined over~$\QQ$, 
the equations that we compute have coefficients in~$\QQ$.
We work this out for the modular curves $X_{\rm ns}^+(p)$ for $p=17$, $19$ and $23$.  

\smallskip
In the remainder of this introduction, we provide some context for  our computational results.  
The curves $X_{\rm ns}(p)$ and $X_{\rm ns}^+(p)$ are defined over~$\QQ$.
Their genera  grow rapidly with~$p$. See~[\BA].  
This may explain why thus far not many computations have 
been done with these curves. 

The curves $X_{\rm ns}(p)$  have no real and hence no rational points.
For $p\le 5$ the genus of  $X_{\rm ns}(p)$ is zero. The curve 
$X_{\rm ns}(7)$ has genus $1$ and, for the record,  is given by the equation $-y^2 = 2x^4 - 14x^3 + 21x^2 + 28x + 7$.
Equations for the genus $4$  curve $X_{\rm ns}(11)$ are given in~[\DGGS]. Using the methods explained in this paper,
 equations for the genus $8$ curve $X_{\rm ns}(13)$ are determined in~[\DMS]. No explicit equations have been computed  for the curves~$X_{\rm ns}(p)$ for primes~$p>13$.

The curves $X_{\rm ns}^+(p)$  are  quotients of $X_{\rm ns}(p)$ by a modular involution. 
The  rational points of the curves $X^+_{\rm ns}(p)$ are relevant  in connection with Serre's Uniformity Conjecture~[\SE]. Indeed, 
after  Mazur's 1978 result ~[\MAI] and the 2010 paper by Bilu, Parent and Rebolledo~[\BPR], the 
conjecture would follow, if  for sufficiently large primes $p$, the only rational points of 
the curves $X^+_{\rm ns}(p)$ are CM-points.

For $p\le 7$ the curves $X^+_{\rm ns}(p)$ have genus zero and  have infinitely many rational points. 
For $p=11$ the genus is~$1$ and there are also infinitely many rational points.
An explicit equation was  computed in 1976 by Ligozat~[\LI].
For $p>11$ the genus exceeds $2$ and hence there are only finitely many rational points.
An equation for the genus $3$ curve $X_{\rm ns}^+(13)$  was computed in 2014 by B.~Baran~[\BB].
In this paper we present equations for $X^+_{\rm ns}(p)$ for the primes $p=17,19$ and~$23$.
Recently Balakrishnan and her coauthors~[\BD] used the Chabauty-Kim method to show that the curve $X_{\rm ns}^+(13)$ 
has precisely seven rational points. All these points are CM-points. For $p>13$ it is at present not known whether or not $X_{\rm ns}^+(p)$ admits any rational points that are not CM. For $p=17,19$ and $23$
a quick computer calculation shows that these curves do not admit any non-CM rational  points that have small coordinates in our models. There may very well not be  any. See sections 6, 7~and~8.

In section 2 we fix our notation and recall some of the basic properties of representations of ${\rm GL}_2(\FF_p)$.
In section 3 we do the same for  the various modular curves that play a role.  In section 4 we determine our trace map  for the principal series and the twisted Steinberg representations. In section 5 we do the same for the cuspidal representations. In section 6 we describe in some detail the actual computations for the curve $X^+_{\rm ns}(17)$. In sections 7 and 8 we present the numerical results
for $X^+_{\rm ns}(19)$ and $X^+_{\rm ns}(23)$.

\beginsection 2. Representations of ${\rm GL}_2(\FF_p)$.

Let $p>2$ be a prime.  In this section we fix notation and recall the  basic properties of the representation theory of the group $G={\rm GL}_2(\FF_p)$, on which our computations are based.

The group $G$ acts on the $p+1$ points of the projective line $\PP_1(\FF_p)$ via linear fractional transformations.
A {\it Borel subgroup} is the stabilizer of a point. It is conjugate to the subgroup $B$ of upper triangular matrices
and has order $p(p-1)^2$.
A {\it split Cartan subgroup} of $G$ is the stabilizer of two points. It is  conjugate to
the subgroup $T$ of diagonal matrices. It has  order $(p-1)^2$ and index $2$ in its normalizer~$N$. 
The group $G$ also acts on the $p^2+1$ points of $\PP_1(\FF_{p^2})$.
A {\it non-split Cartan subgroup} of $G$ is the stabilizer of two   points of $\PP_1(\FF_{p^2})$
that are conjugate over~$\FF_p$. 
Any such group is conjugate to the subgroup $T'$ of matrices that fixes the points $\pm\sqrt{u}$, 
where $u$ denotes a non-square in $\FF_p$. Explicitly, we have 
$$
T'\,=\,\{\pmatrix{a&bu\cr b&a\cr}\in G: \hbox{$a,b\in\FF_p$ with $a^2-ub^2\not=0$}\}.
$$
The group $T'$ is cyclic of order $p^2-1$ and has index $2$ in its normalizer $N'$.

\medskip
In this paper we mostly deal with representations $V$ of $G$ for which the 
subgroup of scalar matrices  $Z$ act trivially. These are representations of $G/Z={\rm PGL}_2(\FF_p)$.  
The complex irreducible representations of ${\rm PGL}_2(\FF_p)$ come in {\it four types}~[\BU, \LAN].
There are two $1$-dimensional representations: the trivial character and a quadratic character $\omega$. Both factor through the determinant. 
There are also two irreducible $p$-dimensional representations. To define them, 
we consider the natural action of ${\rm PGL}_2(\FF_p)$ on the ring $A$ of functions  $\phi:{\bf P}_1(\FF_p)\longrightarrow\CC$
given by $\sigma\phi(P)=\phi(\sigma^{-1}(P))$ for $P\in {\bf P}_1(\FF_p)$.
Since the subspace $\CC$ of constant functions is preserved by this action,
${\rm PGL}_2(\FF_p)$ acts on the $p$-dimensional quotient space $V_{\rm st}=A/\CC$.
This representation is irreducible, has dimension $p$ and is called the {\it Steinberg representation} $V_{st}$.
Its twist by $\omega$ is denoted by $V_{\omega}$.  

The irreducible representations of the third type are the {\it principal series} representations~$V_{\mu}$. These are the inductions of  characters $\mu:B/Z\longrightarrow\CC^*$ 
for which $\mu^2\not=1$. Since the characters $\mu$ are trivial on the unipotent subgroup
$$
U\,=\,\{\pmatrix{1&x\cr0&1\cr}:x\in\FF_p\},
$$ 
they can  be viewed as characters of the cyclic group~$T/Z$.
The representations $V_{\mu}$ have dimension $p+1$. Two representations $V_{\mu}$ and $V_{\mu'}$ are isomorphic if and only if~$\mu'=\mu^{\pm 1}$.  There are $(p-3)/2$ mutually non-isomorphic representations of this type.
The  irreducible representations of the fourth type are the {\it cuspidal} ones. They are associated to 
characters $\theta: T'/Z\longrightarrow\CC^*$ for which  $\theta^2\not=1$. These representations have dimension $p-1$ and are denoted by $V_{\theta}$. Two representations $V_{\theta}$ and $V_{\theta'}$ are isomorphic if and only if $\theta'=\theta^{\pm 1}$. 
There are $(p-1)/2$ mutually non-isomorphic representations of the this type. 
See~[\BU, \LAN] for all this. In section 5 we describe explicit models for the representations $V_{\theta}$.

A character $\mu:T/Z\longrightarrow\CC^*$ is called even or odd, depending on whether it is $1$ on the unique element of order $2$ in $T/Z$ or not. Similarly, a character $\theta:T'/Z\longrightarrow\CC^*$ is called even or odd, depending on whether it is $1$ on the unique element of order 2 in $T'/Z$ or not.  Note that the restriction of the quadratic character $\omega$ to  $T/Z$  is even if and only if its restriction to $T'/Z$ is odd.  This happens if and only if $p\equiv 1\zmod 4$. 

The following proposition gives the dimensions of the   $T$-invariant and $T'$-invariant subspaces $V^T$ and $V^{T'}$ of irreducible representations $V$ of ${\rm PGL}_2(\FF_p)$.

\proclaim Proposition 2.1. Let $V$ be an irreducible complex representation of ${\rm PGL}_2(\FF_p)$ that is not $1$-dimensional.
If $V=V_{st}$, then 
$$
\dim V^T\,=\,2,\quad \dim V^{N}\,=\,1,\quad \hbox{and}\quad \dim V^{T'}\,=\dim V^{N'}\,=\,0.
$$
In all other cases we have
$$
\dim V^T\,=\,\dim V^{T'}\,=\,1,\quad\hbox{and}\quad \dim V^N\,=\,\dim V^{N'}\,\le\,1.
$$
Moreover,  we have
$$\dim V^N=\dim V^{N'}=1,\quad\hbox{if and only if}\quad\cases{V=V_{\mu},& with $\mu$ even,\cr
V=V_{\theta},&  with $\theta$ odd,\cr
V=V_{\omega}, &and $p\equiv 1\zmod 4$.\cr}
$$

\smallskip\noindent{\bf Proof.} We recall 
the remarkable isomorphisms of rational $G$-representations 
$$
\QQ[G/T]\,\cong\,\QQ[G/T']\times V_{\rm st} \times V_{\rm st},\quad\hbox{and}\quad
\QQ[G/N]\,\cong\,\QQ[G/N']\times V_{\rm st},
$$
described by De Smit and Edixhoven in~[\DSED]. 

When $V\not= V_{st}$, the fact that  the vector spaces $V^H$ 
and ${\rm Hom}_G(\QQ[G/H],V)$ are naturally isomorphic for every subgroup $H$ of $G$,
implies that
$\dim V^T=\dim V^{T'}$ and $\dim V^N=\dim V^{N'}$.
To show that $\dim V^T=1$, we observe that $\dim V^T$ is equal to the scalar product
$\langle {\rm Res}_T(\chi_V),1_T\rangle_T$.
Here $\chi_V$ denotes the character of $V$ and $1_T$ is the trivial character on~$T$. 
A standard character computation shows this to be equal to $1$ in all cases.
A similar computation shows that $\langle {\rm Res}_N(\chi_V),1_N\rangle_N$ is $0$ or $1$ 
depending on the parity of the relevant character~$\mu$, $\theta$ or $\omega$.
These computations are particularly straightforward when $V=V_{\mu}$ or~$V_{\omega}$. 
For the cuspidal representations $V=V_{\theta}$,   everything can be computed using the 
description of $V_{\theta}$ as a virtual representation as in~[\BU, \LAN].
Alternatively, one may use the explicit models for $V_{\mu}$ and $V_{\theta}$ given in sections~4~and~5.

For the Steinberg representation $V_{st}$, an explicit calculation shows that
$\dim V_{st}^T=2$ and  $\dim V_{st}^{N}=1$. 
The result by De Smit and Edixhoven implies therefore that
$V_{st}^{T'}$ and  $V_{st}^{N'}$ vanish.

This proves the proposition.

\medskip
In the next sections we construct $T'$-invariant elements in $G$-representations $V$ by applying the
$T'$-trace 
$$
\sum_{t\in T'}t\,=\,\sum_{a,b\in\FF_p,\,a^2-ub^2\not=0}\pmatrix{a&bu\cr b&a\cr}\quad\hbox{in $\ \ \QQ[G]$}
$$
to suitable vectors $v\in V$. 
Since we have the Bruhat decomposition  $G=B\cup BwB$, where
$$
w\,=\,\pmatrix{0&-1\cr1&0\cr},
$$
every non-scalar element in $T'$ can be written as an element in $BwB$. This leads to
the following formula for a  projective version of the $T'$-trace.

\proclaim Proposition 2.2.  The $T'$-trace element $\sum_{M\in T'/Z} M$ of the group ring~$\QQ[{\rm PGL}_2(\FF_p)]$
is given by
$$
{\rm id}\,+\,\sum_{r\in\FF_p}\pmatrix{1&r\cr0&1\cr}w\pmatrix{1&r\cr 0&r^2-u\cr}.
$$

\noindent{\bf Proof.} Representatives in $T'$ of the quotient group $T'/Z$ are the identity matrix
and the matrices  $\pmatrix{r&u\cr 1&r\cr}$  with $r\in \FF_p$. Since  $\pmatrix{r&u\cr 1&r\cr}=
\pmatrix{1&r\cr0&1\cr}w\pmatrix{1&r\cr 0&r^2-u\cr}$, the result follows.

\beginsection 3. Modular curves.

Let $p>2$ be prime and put $G={\rm GL}_2(\FF_p)$. The modular curve $X(p)$ is an algebraic  curve over $\QQ$ that parametrizes elliptic curves with full level $p$ structure. The field of constants of its function field  is the cyclotomic field $\QQ(\zeta_p)$.  The curve $X(p)$ admits a natural morphism to the $j$-line $X(1)$ over~$\QQ$.  The Galois group of $X(p)$ over $X(1)$ is naturally isomorphic to ${\rm GL}_2(\FF_p)/\{\pm{\rm id}\}$.  
Restriction of  automorphisms in ${\rm Gal}(X(p)/X(1))$ to 
the Galois group  of $\QQ(\zeta_p)$ over $\QQ$ coincides with the  determinant map $G/\{\pm {\rm id}\}\longrightarrow\FF_p^*$.

For every subgroup $H$ of $G$ containing $\{\pm {\rm id}\}$ we write $X(p)_H$ for the quotient of $X(p)$ by $H$.
The field of constants of its function field  is the subfield of $\QQ(\zeta_p)$
that is invariant under the subgroup $\det(H)$ of $\FF_p^*$.  We put
$$\Gamma_H\,=\,\{A\in{\rm SL}_2(\ZZ): A\zmod{p}\in H\}.
$$ 
Then the non-cuspidal complex points of any base change of $X(p)_H$ from its field of constants to $\CC$, form the Riemann surface $\Gamma_H\backslash{\bf H}$. 

Taking for $H$ the subgroup $Z$ of scalar matrices of $G$, we obtain the curve $X(p)_Z$. We denote it by $X(p)'$.
Its field of constants  is the quadratic subfield of $\QQ(\zeta_p)$. This is $\QQ(\sqrt{p})$ or $\QQ(\sqrt{-p})$
depending on whether $p\equiv 1$ or $3\zmod 4$.  Since $Z\cap {\rm SL}_2(\FF_p)=\{\pm{\rm id}\}$, the base change
of $X(p)'$ from  $\QQ(\sqrt{\pm p})$ to $\QQ(\zeta_p)$ is the curve $X(p)$.  The curves $X(p)_T$ and $X(p)_N$ 
associated to the split Cartan subgroup $T$ and its normalizer $N$  and the curves $X(p)_{T'}$ and $X(p)_{N'} $   associated to the non-split Cartan subgroup $T'$ and its normalizer $N'$ are quotients of $X(p)'$. These are the  curves $X_{\rm s}(p)$, $X_{\rm s}^+(p)$,  $X_{\rm ns}(p)$ and $X_{\rm ns}^+(p)$ respectively, that were mentioned in the introduction.
Since the determinant maps from the subgroups $T,N,T'$ and $N'$ to $ \FF_p^*$ are all surjective, the curves are all  defined over~$\QQ$, in the sense that their fields of constants are equal to~$\QQ$.

The  group $G={\rm GL}_2(\FF_p)$ acts naturally and linearly on the $\QQ$-vector space $\Omega^1(X(p))$  of K\"ahler differentials.
Therefore its quotient  $G/Z={\rm PGL}_2(\FF_p)$ acts on the $\QQ$-vector space~$\Omega^1(X(p))^Z$ of $Z$-invariants.
On the other hand,   the index $2$ subgroup  ${\rm PSL}_2(\FF_p)$ of ${\rm PGL}_2(\FF_p)$ is isomorphic to the quotient group ${\rm SL}_2(\ZZ)/\Gamma_Z$. Therefore it acts naturally on the complex vector space $S_2(\Gamma_Z)$ of weight $2$ cusp forms for the congruence subgroup $\Gamma_Z$. The two actions are related by the fact that 
 $\Omega^1(X(p)')\otimes_{\QQ}\CC$ is  isomorphic to the induction from ${\rm PSL}_2(\FF_p)$ to ${\rm PGL}_2(\FF_p)$ of~$S_2(\Gamma_Z)$. See~[\BB,~p.279].  So we can write 
 $\Omega^1(X(p)')\otimes_{\QQ}\CC=S_2(\Gamma_Z)+[R]S_2(\Gamma_Z)$ for some 
 fixed respresentative  $R$ of the non-trivial coset of the normal subgroup ${\rm PSL}_2(\FF_p)$ of ${\rm PGL}_2(\FF_p)$.
 Following [\BB], we call  the first coordinate $f_1$ of an element $f_1+[R]f_2$ of  $S_2(\Gamma_Z)+[R]S_2(\Gamma_Z)$,
 its {\it classical coordinate}.

\proclaim Proposition 3.1.  Let $H$ be a subgroup of ${\rm GL}_2(\FF_p)$ containing $Z$.
\item{(a)} The natural maps
 $$
\Omega^1(X(p)_H)\,\mathop{\longrightarrow}\limits^{\cong}\, \Omega^1(X(p))^H\,=\,\Omega^1(X(p)')^{H'},
 $$
 are isomorphisms. Here $H'$ denotes the subgroup $H/Z$ of ${\rm PGL}_2(\FF_p)$. 
 \item{(b)}  If $H$ has the property that $\det(H)=\FF_p^*$, then projection on the classical coordinate
 induces an isomorphism
 $$
\Omega^1(X(p)_H)\otimes_{\QQ}\CC\,\mathop{\longrightarrow}\limits^{\cong}\,S_2(\Gamma_H)
$$
of ${\rm SL}_2(\FF_p)$-representations.
\item{(c)}  Let $H$  be the standard Borel subgroup $B$. It acts on 
$\Omega^1(X_U)\otimes_{\QQ}\CC$ and for any character $\mu$ of $B$,
projection on the classical coordinate induces an isomorphism
$$
(\Omega^1(X_{ZU})\otimes_{\QQ}\CC)(\mu)\mathop{\longrightarrow}\limits^{\cong}\,S_2(\Gamma_1(p),\mu^2).
$$
Here the left hand side denotes the subspace of $\Omega^1(X_{ZU})\otimes_{\QQ}\CC$ on which $B$ acts via the character~$\mu$.
The right hand side is the subspace of $S_2(\Gamma_1(p))$ on which the diamond operators act through the character $\mu^2$.
 
\noindent{\bf Proof.} Part (a) is well known. Part (b) follows from the fact that $H$-invariant elements in $\Omega^1(X(p)')\otimes_{\QQ}\CC=S_2(\Gamma_Z)+[R]S_2(\Gamma_Z)$ are determined by their classical coordinates.  Indeed, we may choose the representative $R$  inside $H$. Then the two coordinates must be equal.

(c) The two coordinates of an element of $\Omega^1(X_{ZU})$ are cusp forms  in $S_2(\Gamma_1(p))$.
The diamond operators in $\Gamma_0(p)/\pm\Gamma_1(p)$ are congruent modulo $p$ to  matrices of the form
$$\pmatrix{a&0\cr0&a^{-1}\cr},\qquad\hbox{with $a\in\FF_p^*$}.
$$
It follows that, if  $b\in B$ acts as multiplication by  $\mu(b)$ on  an element of 
$\Omega^1(X_{ZU})$, the two coordinates are in $S_2(\Gamma_1(p),\mu^2)$. 
The second coordinate is determined by the classical one.  Indeed, we can choose $R\in B$ and then the second coordinate is
equal to the first multiplied by $\mu^{-1}(R)$.

 This proves the proposition
 \medskip

Of special interest is the standard split Cartan subgroup $T$ of $G$. Since the subgroup $\Gamma_T$ of ${\rm SL}_2(\RR)$ is conjugate to $\Gamma_0(p^2)$, there is a natural Hecke compatible isomorphism  $S_2(\Gamma_0(p^2))\longrightarrow S_2(\Gamma_T)$. In terms of $q$-expansions at infinity, the  isomorphism is given by
$$
\sum_{n\ge 1}a_nq^{pn}\,\,\mapsto\,\,\sum_{n\ge 1}a_nq^{n},
$$
where $q$ denotes $\exp(2\pi i\tau/p)$ with $\tau\in{\bf H}$.  Since, the Fourier coefficients of $\Gamma_0(p^2)$-invariant normalized eigenforms are totally real algebraic integers, so are those of $T$-invariant  normalized eigenforms.

We denote the subspace of {\it newforms}  of  $S_2(\Gamma_0(p^2))$  by $S_2(\Gamma_0(p^2))^{\rm new}$.
Abusing notation somewhat, we denote the corresponding subspace of $S_2(\Gamma_T)$  by $S_2(\Gamma_T)^{\rm new}$.  
Note however,  that all forms in $S_2(\Gamma_T)$  are of level $p$. See~[\BB,~(3.4)]. 
By Prop.~3.1~(b) applied to $H=T$, we may identify $S_2(\Gamma_T)$ with the subspace of $T$-invariants of $\Omega^1(X(p))\otimes_{\QQ}\CC$.
By $V_f$ we denote the $\QQ[G]$-subrepresentation of $\Omega^1(X(p))$ generated by a 
normalized eigenform $f$ in $S_2(\Gamma_T)^{\rm new}$. 
It is a vector space over the number field $K_f$ generated by the Fourier coefficients of~$f$.

\proclaim Proposition 3.2.   Let $f$ be a normalized eigenform in $S_2(\Gamma_T)^{\rm new}$. 
Then the subgroup $Z$ of scalar matrices acts trivially on the $\CC[G]$-module $V_f\otimes_{K_f}\CC$.
Moreover,  $V_f\otimes_{K_f}\CC$ is an irreducible  representation of dimension $\not=1$, which is not isomorphic
to the Steinberg representation.

\smallskip\noindent{\bf Proof.}   See~[\BB,~Prop.3.6].  Let $V$ be an irreducible constituent of $V_f\otimes_{K_f}\CC$.
By semi-simplicity we have that $V^T\not=0$. The $G$-action and the Hecke action on 
$\Omega^1(X(p))$ commute. Therefore,  for a prime number $l\not=p$ the Hecke operator $T_l$ acts on $V$  as multiplication 
by the Fourier coefficient~$a_l$. Then it also acts this way on the subspace $V^T$ of $S_2(\Gamma_T)$.
Since $f$  corresponds to a newform in $S_2(\Gamma_0(p^2))$, strong multiplicity one implies that 
$V^T$ is the $1$-dimensional complex vector space generated by~$f$. It follows that $f\in V$, so that $V$ is equal to the 
irreducible representation $V_f\otimes_{K_f}\CC$.   The group $Z$ acts trivially on $V_f$, since it is contained in~$T$.

If $V_f$ had dimension $1$, it  would be invariant under ${\rm PSL}_2(\FF_p)$.
Since the quotient of $X(p)'$ by ${\rm PSL}_2(\FF_p)$ is a  genus $0$ curve over $\QQ(\sqrt{\pm p})$,  
the ${\rm SL}_2(\FF_p)$-invariants of $\Omega^1(X(p)')$ are zero and $V_f$ must be zero as well. Contradiction.
Since the subspace of $T$-invariant elements of $V_f$ has dimension $1$, Proposition 2.1 implies that $V_f$ cannot be the Steinberg representation either.

This proves the proposition.

\beginsection 4. Principal series and twisted Steinberg representations.

In this section we explain how to find elements that are invariant under a non-split Cartan subgroup in the  representations $V_f$
generated by a normalized eigenform $f$, in case  $V_f$ is a principal series 
or twisted Steinberg representation of $G={\rm GL}_2(\FF_p)$ on which the center $Z$ acts trivially.

Let $p>2$ be a prime and let $B$ be the standard Borel subgroup of $G$.
The $1$-dimensional characters of the group $B$ that are trivial on the center $Z$ form a cyclic group of order~$p-1$.  
Given such a character $\mu$, we write $\QQ(\mu)$ for the number field generated by the values of $\mu$.
An explicit model for the  
induced representation ${\rm Ind}_B^G(\mu)$ of $G$ is 
$$
\{\phi:G\longrightarrow \QQ(\mu): \hbox{$\phi(gb)=\mu^{-1}(b)\phi(g)$ for all $g\in G$ and $b\in B$}\}.
$$
The group $G$ acts on this $\QQ(\mu)$-vector space as follows
$$
(\sigma \phi)(x)\,=\,\phi(\sigma^{-1}x),\qquad\hbox{for $\sigma,x\in G$ and  $\phi \in {\rm Ind}_B^G(\mu)$.}
$$
A basis of ${\rm Ind}_B^G(\mu)$ is given by $e_r$ with $r\in\PP_1(\FF_p)=\FF_p\cup\{\infty\}$, where
$e_r$ is the function on $G$ that is equal to $\mu^{-1}$ on the $B$-coset $\{\sigma\in G:\sigma(\infty)=r\}$ and zero elsewhere.
For every $r\in\FF_p$ the $G$-action on $e_r$ can easily be computed: for every $k\in\FF_p$ we have
$$
\pmatrix{1&k\cr0&1\cr}e_r=e_{r+k},\quad\hbox{for $r\in\FF_p$, while} \quad\pmatrix{1&k\cr0&1\cr}e_{\infty}=e_{\infty}.\eqno{(1)}
$$
For every $a\in\FF_p^*$ we have
$$
\pmatrix{a&0\cr0&1\cr}e_r=\mu\pmatrix{1&0\cr0&a\cr}e_{ar}\quad\hbox{for $r\in\FF_p$, while} 
\quad\pmatrix{a&0\cr0&1\cr}e_{\infty}=\mu\pmatrix{a&0\cr0&1\cr}e_{\infty}.\eqno{(2)}
$$
The action of the matrix $w=\pmatrix{0&1\cr-1&0\cr}$ is given by 
$$
we_r=\mu\pmatrix{1/r&0\cr0&r\cr}e_{-1/r}\quad\hbox{for $r\in\FF_p^*$},\eqno{(3)}
$$
while $w$ switches $e_0$ and $e_{\infty}$. Since $G=B\cup BwB$, these formulas determine the action of~$G$.

If $\mu^2\not=1$,  we  recover the irreducible complex representation  $V_{\mu}$ of section 2 as ${\rm Ind}_B^G(\mu)\otimes_{\QQ(\mu)}\CC$.   The values of the character of $V_{\mu}$ generate the maximal real subfield $\QQ(\mu)^+$ of the cyclotomic field~$\QQ(\mu)$.
Since the subspace of $T$-invariants is $1$-dimensional, it follows from~[\WAL,~Lemma 1.1] that the representation $V_{\mu}$ itself can actually be defined over $\QQ(\mu)^+$. We do not make use of this.

If $\mu^2=1$, we have that $\mu$ is either $1$ or $\omega$, so that $\QQ(\mu)=\QQ$. 
In this case,  the subspace $L$  of constant functions of ${\rm Ind}_B^G(\mu)$, generated by $e_{\infty}+\sum_{r\in\FF_p}e_r$, is preserved by~$G$ and the representation $({\rm Ind}_B^G(\mu)/L)\otimes_{\QQ}\CC$ is irreducible. 
In fact, we recover the complex Steinberg representation $V_{\rm st}$ and its quadratic twist $V_{\omega}$. See~[\BU].

We let $Z, B,T,T',N, N'$ and $U$ be the subgroups of $G$ defined in section 2. 
It is convenient to view $\mu$ as a character of $\FF_p^*$. For this reason we put
$$
\mu(r)\,=\,\mu\pmatrix{r&0\cr0&1\cr},\qquad\hbox{for $r\in\FF_p^*$}.
$$

\proclaim Proposition 4.1. Let $\mu:B/Z\longrightarrow \CC^*$ be a  character satisfying $\mu^2\not=1$
and let $V_{\mu}$ be the principal series representation associated to~$\mu$.

\item{(a)} The subspace of $V_{\mu}$ of $U$-invariants  has dimension $2$ and is generated by $e_{\infty}$  and
by $\sum_{r\in\FF_p}e_r$. The subgroup $B$ acts  via $\mu$ on the line generated by $e_{\infty} $
and via $\mu^{-1}$ on the line generated by~$\sum_{r\in\FF_p}e_r$.  
\item{(b)} The subspace of $T$-invariants is generated by 
$$\sum_{r\in\FF_p^*}\mu(r)e_r.
$$
It is invariant under the action of the normalizer $N$ if and only if $\mu$ is an even
character of~$B/ZU=T/Z$.
\item{(c)} The subspace of $T'$-invariants is generated by 
$$
e_{\infty} +\sum_{r\in\FF_p}\mu^{-1}(r^2-u)e_r.
$$
It is invariant under the action of the normalizer $N'$ if and only if $\mu$ is even.

\smallskip\noindent {\bf Proof.} Parts (a) and (b)  easily follow from the formulas given above. 
The computations are easy and left to the reader. 
By Proposition 2.1, the subspaces of $T$-invariants and of $T'$-invariants have dimension~1.
The element listed in (c) is  the $T'$-trace of Proposition 2.2 applied to $e_{\infty}$. 

This proves the proposition.

\medskip
For the character $\mu=\omega$, the result is similar:

\proclaim Proposition 4.2. Let  $\omega$ be the quadratic character of $G$
and let $V_{\omega}$ be the  twisted Steinberg representation.
\item{(a)} The subspace of $V_{\omega}$ of $U$-invariants  has dimension $1$ and 
is generated by $e_{\infty}$. 
The subgroup $B$ acts   on it via $\omega$.  
\item{(b)} The subspace of $T$-invariants is generated by 
$$\sum_{r\in\FF_p^*}\omega(r)e_r.
$$
It is invariant under the action of the normalizer $N$ if and only if  $p\equiv 1\zmod 4$.
\item{(c)} The subspace of $T'$-invariants is generated by 
$$
e_{\infty} +\sum_{r\in\FF_p}\omega(r^2-u)e_r.
$$
It is invariant under the action of the normalizer $N'$ if and only $p\equiv 1\zmod 4$.

\smallskip\noindent {\bf Proof.} Note that in $V_{\omega}$ we have the relation $e_{\infty}=-\sum_{r\in\FF_p}e_r$.
The proof is similar to the proof of Proposition 4.1.

\medskip
The relation with $q$-expansions of eigenforms is as follows.  Suppose that $f$ is a normalized eigenform  in the space $S_2(\Gamma_T)^{\rm new}$. Since $f$ is $\Gamma_T$-invariant, Proposition 3.1 (b) implies that 
we can identify $S_2(\Gamma_T)^{\rm new}$ with a subspace of 
$\Omega^1(X(p))\otimes_{\QQ} \CC$. 
By Proposition~3.2, the newform $f$ generates an absolutely irreducible  $G$-representation~$V_f$,  defined over 
the number field $K_f$ generated by the Fourier coefficients of $f$. 
Note that this implies that $K_f$ contains the field $\QQ(\mu)^+$ of character values.

Suppose that $V_f$ is a principal series or twisted Steinberg representation. In other words, we have an isomorphism
$$
V_{\mu}\,\cong\,V_{f}\otimes_{K_f}\CC,\qquad\hbox{for  some non-trivial character $\mu:B/Z\longrightarrow\CC^*$.}
$$
By Propositions 4.1 (a) and 4.2 (a), the representation $V_{\mu}$ admits a unique $1$-dimensional $U$-invariant subspace $W$
on which the Borel  subgroup $B$ acts via $\mu$.  It is generated by the element~$e_{\infty}$.
Proposition 3.1 (c) implies then that in $V_f$, there is a unique element whose  classical coordinate is a  $\Gamma_1(p)$-invariant normalized eigenform~$h$ on which $\Gamma_0(p)$ acts via the character~$\mu^2$.  In the twisted Steinberg case, we have $\mu=\omega$ and hence $\mu^2=1$. In this case $h$ is a $\Gamma_0(p)$-invariant normalized eigenform.   

Any $G$-equivariant linear map  $V_{\mu}\longrightarrow V_f\otimes_{K_f}\CC$, must map $e_{\infty}$ into the $1$-dimensional space generated by $h$.
Schur's Lemma  implies that for each $c\in\CC^*$ there is a unique $G$-equivariant isomorphism
$$
j_c:V_{\mu}\mathop{\longrightarrow}\limits^{\cong} V_f\otimes_{K_f}\CC,
$$
for which $j_c(e_{\infty})=c h$.

Let $q=e^{{2\pi i\tau}\over p}$. Since  $h$ is $\Gamma_1(p)$-invariant,  its Fourier expansion is  of the form
$$
h\,=\,\sum_{n\ge 1} a_nq^{pn}.
$$  
Note that there is also a unique element  in  $V_f\otimes_{K_f}\CC$ whose classical coordinate is the `complex conjugate' normalized eigenform~$\overline{h}=\sum_{n\ge 1} \overline{a_n}q^{pn}\in S_2(\Gamma_1(p),\mu^{-2})$. The isomorphism $j_c$ maps the element $-\sum_{r\in\FF_p}e_r$ to
a multiple of $\overline{h}$.

The following proposition relates the Fourier expansion of $f$ to the one of $h$.

\proclaim Proposition 4.3.  Let $\mu\not=1$ and let $f$ and $h$ be the normalized eigenforms described  above. Put $\zeta_p=e^{{2\pi i}\over p}$.
\item{(a)}
Then the  $q$-expansion of $f$  is given by
$$
f\,=\,\sum_{n\ge 1}\mu(n)a_nq^n,
$$
with the convention that $\mu(n)=0$, whenever $n$ is divisible by~$p$. 
\item{(b)}  The eigenform $h$ is in the $\QQ(\mu)[G]$-span of ${{\tau(\mu)\tau(\mu^2)}\over{a_p}}f$. Here $\tau(\mu)$ and $\tau(\mu^2)$ denote the  Gaussian sums
$\sum_{x\in\FF_p}\mu(x)\zeta_p^x$ and $\sum_{x\in\FF_p}\mu^2(x)\zeta_p^x$ respectively.
When $\mu=\omega$ we have $\mu^2=1$ and we put $\tau(\mu^2)=-1$.

\smallskip\noindent{\bf Proof.}   By Prop.~4.1 (b), the subspace of $T$-invariant elements of $V_{\mu}$
is the $1$-dimensional subspace generated by  $\sum_{r\in\FF_p^*}\mu(r)e_r$. The isomorphism $j_c$ introduced above, maps it to a $\Gamma_T$-invariant eigenform  in $V_{f}\otimes_{K_f}\CC$. For a suitable choice of $c$ we obtain $f$ itself.

We compute $j_c(\sum_{r\in\FF_p^*}\mu(r)e_r)$. The formulas (1), (2) and (3) given above, imply  that 
$$
e_r\,=\,\pmatrix{1&r\cr0&1\cr}w e_{\infty},\qquad\hbox{for $r\in\FF_p$.}
$$
It follows that
$$
j_c(e_r)\,=\,c\pmatrix{1&r\cr0&1\cr}wh.
$$
By Atkin-Li~[\AL], the modular involution $w_p$ transforms $h$ into the `complex conjugate' form $\overline{h}$ multiplied by the so-called pseudo-eigenvalue~$\epsilon$, which is a complex number of absolute value~$1$. To be precise, $\epsilon$ is equal to $\tau(\mu^2)/a_p$ and we have
$$
{1\over{p\tau^2}}h(-{1\over{p\tau}})\,=\,\epsilon\overline{h}(\tau),\qquad\hbox{for $\tau\in{\bf H}$.}
$$
This  implies 
that $wh$ is the element   of $V_{\mu}$ whose classical coordinate
is equal to the Fourier series
$$
wh(\tau)\,=\, {{\epsilon}\over p}\overline{h}(\tau/p)\,=\,{{\epsilon}\over p}\sum_{n\ge 1} \overline{a_n}q^{n}.
$$
It follows that for $r\in\FF_p$ we have
$$
j_c(e_r)\,=\, c\pmatrix{1&r\cr0&1\cr}wh(\tau)\,=\, c{{\epsilon}\over p}\sum_{n\ge 1} \overline{a_n}\zeta_p^{rn}q^{n}.
$$
Therefore the classical coordinate of $j_c(\sum_{r\in\FF_p^*}\mu(r)e_r)$  is
$$
c{{\epsilon}\over p}\sum_{n\ge 1}\sum_{r\not=0} \overline{a_n}\mu(r)\zeta_p^{nr}q^n\,=\,c
{{\epsilon\tau(\mu)}\over p}\sum_{n\ge 1}\mu^{-1}(n)\overline{a_n}q^n
\,=\,c
{{\epsilon\tau(\mu)}\over p}\sum_{n\ge 1}\mu(n)a_nq^n.
$$
The last equality follows from the fact that $\mu^{-1}(n)\overline{a_n}$ is real and hence equal to $\mu(n)a_n$ for all $n\in\ZZ$.

Since $f$ is a {\it normalized} eigenform, part (a) follows. 
When we choose $c=(\epsilon\tau(\mu)/p)^{-1}$,
we have that $j_c(\sum_{r\in\FF_p^*}\mu(r)e_r)=f$. In particular,
$f$ is in the $\QQ(\mu)[G]$-span of $j_c(e_{\infty})=ch$.  Since $V_f$ is irreducible, this is the same as saying that $h$ is in the $\QQ(\mu)[G]$-span of $\epsilon\tau(\mu)f$.

This proves the proposition.

\medskip
We now turn to the computation of the Fourier series of the $T'$-invariant eigenform in~$V_{\mu}$.
See also~[\BAO].
Recall that $u\in\FF_p^*$ is a fixed non-square. We put
$$
\lambda_n=\sum_{r\in\FF_p}\mu^{-1}(r^2-u)\zeta_p^{rn},\qquad\hbox{for $n\in\ZZ$.} 
$$

\proclaim Proposition 4.4.   Let $f\in S_2(\Gamma_T)^{\rm new}$  be the $T$-invariant eigenform discussed above
and let $h=\sum_{n\ge 1}a_nq^{pn}$ be the corresponding $\Gamma_1(p)$-invariant eigenform. Then the element of $V_f$ with classical coordinate equal to 
$$
{1\over{\tau(\mu)}}\left({p\over{\tau(\mu^{2})}}\sum_{n,\,p|n}a_nq^n\,+\,\sum_{n\ge 1}\lambda_n\overline{a_n}q^n\right)
$$
is a generator for the subspace of $T'$-invariant forms.
Moreover, it is in the $\QQ(\mu)[G]$-span of the $\Gamma_T$-invariant eigenform $f$.

\smallskip\noindent{\bf Proof.}   Propositions 4.1 (c) and 4.2 (c) give an explicit generators of the $1$-dimensional subspace of $T'$-invariants of $V_{\mu}$.  We apply the ismorphism $j_c$ with $c=p/\epsilon\tau(\mu)$ as we did above. Since
$$
e_r\,=\,\pmatrix{1&r\cr0&1\cr}w e_{\infty},\qquad\hbox{for $r\in\FF_p$,}
$$
we get
$$
{p\over{\epsilon\tau(\mu)}}\left(\sum_{n\ge 1}a_nq^{pn}+{{\epsilon}\over p}\sum_{r\in\FF_p}\mu^{-1}(r^2-u)\sum_{n\ge 1}\overline{a_n}\zeta_p^{rn}q^n\right).
$$
Since $a_{pn}=a_pa_n$ for every $n\ge 1$, this is equal to
$$
{p\over{\tau(\mu)\tau(\mu^2)}}\sum_{n,\,p|n}a_nq^n\,+\,{1\over{\tau(\mu)}}\sum_{n\ge 1}\left(\sum_{r\in\FF_p}\mu^{-1}(r^2-u)\zeta_p^{rn}\right) \overline{a_n}q^n,
$$
which is easily seen to give the result.
By Proposition 4.3 (b) the result is contained  in the $\QQ(\mu)[G]$-span of $f$.
This proves the proposition.

\medskip
The numbers $\lambda_n=\sum_{r\in\FF_p}\mu^{-1}(r^2-u)\zeta_p^{rn}$  are so-called {\it Sali\'e sums}. They are related to Kloosterman sums. See~[\KCON] and the references therein.

\beginsection 5. Cuspidal Representations.

Let $p>2$ be prime,  let $G={\rm GL}_2(\FF_p)$ and let $T\subset G$ be the standard split Cartan subgroup of diagonal matrices.
In this section we consider normalized $\Gamma_T$-invariant weight $2$ eigenforms $f$, that generate representations $V_f\subset \Omega^1(X(p))$  that are cuspidal. We  explain how to find  elements in $V_f$ that are invariant under a non-split Cartan subgroup of $G$.

Let $u\in\FF_p^*$ be a non-square, let $T'$ denote the non-split torus in $G$ introduced in section 2 
and let $\theta: T'\longrightarrow \QQ(\theta)^*$ be a character that is trivial on the subgroup $Z$ of scalar matrices.
We have $\theta^{p+1}=1$ and assume that $\theta^2\not=1$.    By $\QQ(\theta)$  we denote the field generated by the image of $\theta$.
Our model $V_{\theta}$ for the cuspidal representation associated to $\theta$ is the quotient of the $\QQ(\theta)$-vector space of functions $\phi:\FF_p\longrightarrow \QQ(\theta)$ by the $1$-dimensional subspace of constant functions.
The standard Borel subgroup $B\subset G$ acts by fractional linear transformations on $\FF_p={\bf P}_1(\FF_p)-\{\infty\}$
and hence on the space of functions $\phi:\FF_p\longrightarrow \QQ(\theta)$: we have
$\sigma \phi(x) = \phi(\sigma^{-1}x)$ for $\sigma\in B$ and any function $\phi$.
Since  $B$ preserves the constant functions, it acts on $V_{\theta}$.

It is easy to see that $V_{\theta}$ is an  irreducible $(p-1)$-dimensional representation of  $B$, on which the  scalar matrices act trivially.
We turn $V_{\theta}$ into  an irreducible representation of  ${\rm PGL}_2(\FF_p)$.  Let 
$$
w=\pmatrix{0&1\cr-1&0\cr}
$$
be the usual involution. Since $G=B\cup BwB$, it suffices to describe the action of~$w$.
It is given by
$$
w\,\phi =-{1\over p}\sum_{y\in \FF_{p^2}^*}\theta(y)\pmatrix{{\rm N}(y)&{\rm Tr}(y)\cr0&1\cr}\phi,\qquad\hbox{for all $\phi\in V_{\theta}$.}
$$
Here $\FF_{p^2}$ denotes $T'\cup\{\pmatrix{0&0\cr0&0\cr}\}$. It is a subfield of the ring of $2\times 2$ matrices over~$\FF_p$.
Let ${\rm N}$ and ${\rm Tr}$ denote the norm and trace maps from $\FF_{p^2}$ to $\FF_p$ respectively.

 \medskip
Proving that the formula for the action of $w$ gives rise to a well defined action of $G$ on $V_{\theta}$ is straightforward but somewhat cumbersome.  Alternatively, one can relate $V_{\theta}$ to the  representation space described  by 
Bump~[\BU,~4.1].   Let $\zeta_p$ denote a fixed $p$-th root of unity. To every $\phi\in V_{\theta}$ we
associate the function $\tilde{\phi}:\FF_{p^2}^*\longrightarrow \QQ(\theta)$ given by $\tilde{\phi}(y)=\theta^{-1}(y)\sum_{r\in\FF_p}\phi(r)\zeta_p^{r{\rm N}(y)}$. This gives an isomorphism of  $V_{\theta}\otimes_{\QQ(\theta)}\CC$ with Bump's model.
Our model has the advantage  that it can be defined over $\QQ(\theta)$, rather than over a field that contains the $p$-th roots of unity.
The character values of $V_{\theta}$ generate the maximal real subfield $\QQ(\theta)^+$ of $\QQ(\theta)$. As in the principal series case, it follows
from [\WAL,~Lemma 1.1] that $V_{\theta}$ can actually be defined over $\QQ(\theta)^+$. We do not make use of this.

Let $e_0:\FF_p\longrightarrow\QQ(\theta)$ be the characteristic function of~$0$.
For $r\in\FF_p$, let  $e_r=\pmatrix{1&r\cr0&1\cr}e_0$. It is the characteristic function of the element~$r\in\FF_p$.
The functions $e_r$, $r\in\FF_p^*$ form a basis for  the $V_{\theta}$. Since $\sum_{r\in\FF_p}e_r$ is the constant function $1$, we have the relation $\sum_{r\in\FF_p}e_r=0$ in $V_{\theta}$.

\proclaim Proposition 5.1.  Let  $\theta: T'/Z\longrightarrow \QQ(\theta)^*$ be a character satisfying $\theta^2\not=1$ and let
$V_{\theta}$ be the cuspidal representation of $G$ associated to the character~$\theta$. Then
\item{(a)} the subspace of $U$-invariants is zero;
\item{(b)} the subspace of $T$-invariants is generated by $e_0$;
it is invariant under the action of the normalizer $N$ of $T$ if and only if $\theta$ is an odd character of the cyclic group $T'/Z$;
\item{(c)} there is an $r\in\FF_p^*$ for which  the element
$$
pe_r-\sum_{m\in\FF_p}\sum_{y\in\FF_{p^2}^*}\theta(y)e_{{{(m+r){\rm N}(y)}\over{m^2-u}}+{\rm Tr}(y)+m}
$$
generates	 the $1$-dimensional subspace of $T'$-invariants. 
 The space of $T'$-invariants is also $N'$-invariant if and only if $\theta$ is odd.

\smallskip\noindent {\bf Proof.}  Part (a) and the first statement of (b) easily follow from the formulas given above.
The statement about the normalizer $N$ can be proved with a short computation~[\BB,~Prop.2.1].
To prove (c),  we combine the formula  for the action of $w$ with  Proposition 2.2.
It follows that the  $T'$-trace is  equal to
$$
{\rm id}-{1\over p}\sum_{y\in\FF_{p^2}^*}\theta(y)\sum_{m\in\FF_p}\pmatrix{{\rm N}(y)&{m{\rm N}(y)}+(m^2-u)({\rm Tr}(y)+m)\cr0&m^2-u\cr}.
$$
Applying it to $pe_r$ gives the  element of part (c). 
Since the elements $e_r$, with $r\in\FF_p$, generate $V_{\theta}$, their $T'$-traces generate 
the $1$-dimensional space of $T'$-invariants. In other words,  the $T'$-trace of at least one  of the elements $e_r$ is not zero and hence generates the subspace of $T'$-invariants.

This proves the proposition.

\medskip
The relation with $q$-expansions of eigenforms is as follows.
 Let $q=e^{{2\pi i\tau}\over p}$, let $\zeta_p=e^{{2\pi i}\over p}$  and let 
$$
f=\sum_{n\ge 1} a_nq^n
$$  
be a normalized weight 2 eigenform that is invariant under~$\Gamma_T$. By Prop.~3.1~(b) we 
may identify $S_2(\Gamma_T)$ with the $T$-invariant elements in $\Omega^1(X(p))\otimes_{\QQ}\CC$.
Then $f$ generates an absolutely irreducible  $G$-representation~$V_f$ that is defined over 
the number field $K_f$ generated by the Fourier coefficients of $f$.
Suppose that $V_f$ is a cuspidal representation. In other words, we have 
$$
V_{\theta}\,\cong\,V_{f}\otimes_{K_f}\CC,\qquad\hbox{for  some character $\theta:T'/Z\longrightarrow\QQ(\theta)^*$ with $\theta^2\not=1$.}
$$
Note that $K_f$ contains the values of the character of $V_{\theta}$. This means that 
$\QQ(\theta)^+$ is a subfield of $K_f$.

By Proposition 5.1 (b), the element $f\in V_f$ corresponds
to the vector $e_0\in V_{\theta}$ or a multiple thereof. More generally, for any $r\in\FF_p$ the elements in $V_f$ with classical coordinate
equal to
$$
f_r=\pmatrix{1&r\cr0&1\cr}f = \sum_{n\ge 1} a_n\zeta_p^{nr}q^{n},
$$
correspond to multiples of~$e_r$.

\proclaim Proposition 5.2.  The elements in $V_f$ with classical coordinate equal to
$$
pf_r-\sum_{m\in\FF_p}\sum_{y\in\FF_{p^2}^*}\theta(y)f_{{{(m+r){\rm N}(y)}\over{m^2-u}}+{\rm Tr}(y)+m},\qquad\hbox{for $r\in\FF_p$},
$$
are all $T'$-invariant. They are all in the $\QQ(\theta)[G]$-span of $f$ and 
at least one of them generates the subspace of $T'$-invariants of $V_f$.

\smallskip\noindent{\bf Proof.}  This follows from the fact that the vectors $e_r$ are in the $\ZZ[G]$-span
of $e_0$ and the fact that the $T'$-trace is an element of $\QQ(\theta)[G]$.

\beginsection 6. Level $17$.

In this section we explain how to compute equations over $\QQ$ for the canonical embedding of the genus $6$ curve~$X_{\rm ns}^+(17)$.  We follow the method in~[\MER].  We exhibit six linearly independent weight $2$ cusp forms that are invariant under the normalizer $N'$ of the standard
non-split Cartan subgroup $T'$. We find these forms inside the six representation spaces $V_f$, generated by six normalized eigenforms $f\in S_2(\Gamma_T(17))^{\rm new}$, that are  invariant under the normalizer $N$ of the standard split Cartan subgroup $T$.  Since the space $S_2(\Gamma_T(17))^{\rm new}$ is closely related to the space $S_2(\Gamma_0(17^2))^{\rm new}$, we start from there.
We can find the Fourier expansions of the normalized eigenforms in~[\ST],  for instance.   

Up to Galois conjugation and twists by the quadratic character  $\omega$, there are four normalized weight $2$ eigenforms   invariant under $\Gamma_0(17^2)$.
Two of these are  twists of normalized eigenforms in $S_2(\Gamma_1(17))$. They give rise to principal series and twisted Steinberg representations. The other two eigenforms generate cuspidal representations.

There is a unique $\Gamma_0(17)$-invariant normalized eigenform $f_0=\sum_na_nq^{17n}$.  Its $17$-th Fourier coefficient $a_{17}$ is equal to~$+1$.
Its quadratic twist $\sum_n\omega(n)a_nq^{17n}$ is a normalized $\Gamma_0(17^2)$-invariant eigenform. Here
we put  $q=e^{{2\pi i\tau}\over{17}}$ for~$\tau\in\HH$. By convention $\omega(n)=0$ whenever $n$ is divisible by $17$. The corresponding $\Gamma_T$-invariant form is $\sum_n\omega(n)a_nq^{n}$. The first few terms of its Fourier expansion are
$$
f_1\,=\,q-q^2 - q^4 + 2q^5 - 4q^7 + 3q^8 - 3q^9 -2q^{10} - 2q^{13} + 4q^{14} - q^{16} +  3q^{18} +\ldots
$$
The irreducible subrepresentation $V_{f_1}$  of $\Omega^1(X(p)')$ is isomorphic to the twisted Steinberg representation.
The form $f_1$ is also  invariant under the normalizer $N$ of $T$ because $17\equiv 1\zmod 4$. See~Prop.~4.2. 

One finds in Stein's tables that the  space $S_2(\Gamma_1(17))$  is the direct product of the $1$-dimensional space of $\Gamma_0(17)$-invariant forms and a $4$-dimensional subspace $W$ spanned by the Galois conjugates of an eigenform $h$ on which the diamond operators act through a character of order~$8$ of~$(\ZZ/17\ZZ)^*$.  
Any such character is of the form $\mu^2$, where $\mu$ has order $16$. Since $\mu$ is then an odd character of $T/Z$, Prop.~4.1 implies that the  normalizer $N$ acts as $-1$ on the $T$-invariants. 

Therefore the twist  by $\mu$ of $h$ as described in section 4,  is a  $\Gamma_0(17^2)$-invariant normalized weight $2$ eigenform, corresponding to a $\Gamma_T$-invariant form that is not $N$-invariant. It plays no role in our computation of the canonical embedding  of~$X_{\rm ns}^+(17)$.

The remaining normalized eigenforms in Stein's tables both generate cuspidal representations.
Put $a={{-1+\sqrt{13}}\over 2}$. Then the modular form
$$\eqalign{
f_2 \,=\,  &\,q -(a+1)q^2+aq^3+(a+2)q^4-(a+1)q^5-3q^6+(a-1)q^7-3q^8-aq^9+\cr&+(a+4)q^{10}-3q^{11}+(a+3)q^{12}-(a+2)q^{13}+(a-2)q^{14}-3q^{15}+(a-1)q^{16} +\ldots\cr}
$$
is the $\Gamma_T$-invariant form associated to a newform in $\Gamma_0(17^2)$. The representation $V_{f_2}$ is cuspidal with respect to some character $\theta$ of order dividing~$18$. Since  the field $K_{f_2}$ generated by the Fourier coefficients, contains $\QQ(\theta)^+$, we actually must have that~$\theta^6=1$.
Figuring out what character $\theta$ of $T'/Z$ is involved, can be done by numerically computing the action of $w$ on $f_2$ 
for every possible $\theta$ in a suitable $\tau\in\HH$ as in Baran's paper~[\BB,~section 6]. It turns out that in this case $\theta$ has order~$6$.
The twist of $f_2$ by $\omega$ is cuspidal with character $\theta\omega$, which has order $3$. By Prop.~5.1 the form $f_2$ is $N$-invariant, while its twist
is anti-invariant.

The fourth normalized eigenform is the $\Gamma_T$-invariant form associated to one of the  $\Gamma_0(17^2)$-invariant eigenforms in Stein's table   with Fourier  coefficients  in $\QQ(\zeta_9)^+$. The first few terms of its Fourier expansion are
$$\eqalign{
f_3\,=\,&q-(b^2+b-2)q^2-(b+1)q^3+bq^4+(b^2+b-4)q^5+(2b^2 + 2b - 3)q^6+bq^7+\cr
+&(b^2 + b - 3)q^8+(b^2 +2b - 2)q^9+(2b^2 + b - 6)q^{10}
-(2b^2 - 2)q^{11}-(b^2 +b)q^{12}+\ldots \cr}
$$
Here $b=\zeta_9^{}+\zeta_9^{-1}$. It is a zero  of $x^3-3x+1$. 
The representation $V_{f_3}$ is cuspidal with character $\theta$ of order $18$. The twist by $\omega$ is cuspidal with character $\theta\omega$ of order $9$. By Prop.~5.1 the form $f_3$ is $N$-invariant, while its twist is anti-invariant.

At this point we have six $T$-invariant eigenforms: $f_1$, $f_2$ and its Galois conjugate and $f_3$ with its two Galois conjugates.
To $f_1$ we apply the $T'$-trace fomula in Proposition 4.4. This gives us a $T'$-invariant form $g_1$ with Fourier coefficients 
in $\QQ(\zeta_{17})^+$. Applying  the formula of Proposition 5.2 to
$f_2$ and its conjugate over $K_{f_2}=\QQ(\sqrt{13})$, we obtain the $T'$-invariant form $f_2'$ and its conjugate.
Their Fourier coefficients are in  $K_{f_2}(\zeta_{17})^+$.  We put $g_2={\rm Tr}(f_2')$  and $g_3={\rm Tr}(\sqrt{13}f_2')$.
Here ${\rm Tr}$ denotes the trace map from $\QQ(\zeta_{17})^+(\sqrt{13})$ to $\QQ(\zeta_{17})^+$.
Then $g_2$ and $g_3$ are $T'$-invariant forms with Fourier coefficients in $\QQ(\zeta_{17})^+$.
Similarly, we apply the $T'$-trace map given in Proposition 5.2  to $f_3$ and its conjugates over 
$K_{f_3}=\QQ(\zeta_9^{})^+$ and obtain the $T'$-invariant form $f_3'$. Its Fourier coefficients are in
$K_{f_3}(\zeta_{17})^+$.  For $i=1,2,3$, we put $g_{3+i}={\rm Tr}(e_if_3')$, where 
$e_1,e_2,e_3$ denotes the basis of $K_{f_3}(\zeta_{17})^+$ over $\QQ(\zeta_{17})^+$
given by  $1,\alpha,\alpha^2$, where $\alpha$ is a zero of  the defining polynomial $x^3+3x^2-3$ used in Stein's table.
Then $g_4$, $g_5$ and $g_6$ are $T'$-invariant forms with Fourier coefficients in $\QQ(\zeta_{17})^+$.

We list the first few Fourier coefficients of the $T'$-invariant forms $g_1,\ldots,g_6$. By an $8$-tuple $[x_1,\ldots,x_8]\in\ZZ^8$ we denote the element
$\sum_{j=1}^8 x_j(\zeta_{17}^j+\zeta_{17}^{-j})$. For every $i$ we have divided the coefficients of $g_i$ by a common divisor in $\ZZ$.

$$\eqalign{
g_1&\!=[7, 1, 2, 5, 4, 5, 4, 6]q-[6, 7, 4, 1, 5, 2, 4, 5]q^2+[-5, 6, 4, 7, 2, 4, 5, 1]q^4\ldots\cr
g_2&\!=[4, 16, 2, -4, -2, 8, -8, 18]q+[9, 2, -4, 8, 4, 1, -1, -2]q^2-[4, -1, 2, -4, -2, 8, 9, 1]q^3...\cr
g_3&\!=[9, 2, -4, 8, 4, 1, -1, -2]q^2-[4, -1, 2, -4, -2, 8, 9, 1]q^3-[-2, 9, -1, 2, 1, -4, 4, 8]q^4\ldots\cr
g_4&\!=[8, 8, -2, 4, 5, -2, -1, -3]q-[3, 2, -1, 2, -2, 7, -4, 10]q^2-[12, 9, 12, 6, 18, 12, 9, 24]q^3.\rlap.\cr
g_5&\!=-[4, 4, 8, 6, 3, 4, 2, 3]q+[1, 4, -1, 4, -2, 4, -1, 8]q^2+[2, 5, 10, 1, 12, 10, 2, 9]q^3\ldots\cr
g_6&\!=[10, 10, 9, 12, 5, 2, 1, 2]q-[5, 12, 0, 12, 0, 16, 1, 22]q^2-[8, 10, 22, 4, 32, 22, 9, 29]q^3\ldots\cr}
$$
By~[\SD] the canonical embedding of a genus $6$ curve is typically cut out by six quadrics. See also~[\DOSE, Thm.~1.1] and~[\MER].
We compute six quadrics that vanish on the canonically embedded curve $X_{\rm ns}^+(17)$ and then use MAGMA to check that 
the intersection of the quadrics is a curve of genus $6$. Then we know that the quadrics are indeed equations for~$X_{\rm ns}^+(17)$.

To do this, we compute Fourier series of the $21$ products $g_ig_j$ with $1\le i\le j\le 6$.
Even though the Fourier coefficients of the forms $g_i$ are in $\QQ(\zeta_{17})^+$ and are usually not rational,  the corresponding K\"ahler differentials are rational. This is explained by the fact that the cusps of  $X_{\rm ns}^+(17)$ are not rational, but conjugate over $\QQ(\zeta_{17})^+$.
Since the curve $X_{\rm ns}^+(17)$ is defined over $\QQ$, we search for quadrics 
$$\sum_{1\le i\le j\le 6}a_{ij}x_ix_j,
$$
with coefficients $a_{ij}$ in~$\QQ$.  From the equation $\sum_{1\le i\le j\le 6}a_{ij}g_ig_j=0$ we obtain infinitely many equations
with coefficients in $\QQ(\zeta_{17})^+$,  one for every term $q^n$ in the Fourier expansion. 
Since the coefficients are in the degree $8$ number field $\QQ(\zeta_{17})^+$,
each equation gives rise to {\it eight} equations with coefficients in~$\ZZ$.
For instance, a consideration of the Fourier coefficients of $q^2$ and $q^3$ gives rise to the following $16$ equations.
Here the columns correspond to the coefficients $a_{ij}$ in lexicographic order.
$$
\def\quad{\hskip1ex\relax}
\matrix{ 
6&0&0&3&-2&5&-3840 &0& 0 &-2 &2 &0 &0 &0 &0 &15 &2 &7 &3 &-3 &10\cr
3&0&0&3&1&1&10620&0&6&-2&4&0&0&0&0&18&2&8&4&-5&14\cr
4&-2&0&-1&0&-1&-5256&0&-4&0&2&0&0&0&0&7&0&4&6&-9&17\cr
5&2&0&1&0&1&2820&0&-14&0&-8&0&0&0&0&20&0&14&6&-9&24\cr
3&0&0&0&1&-2&-3948&0&12&10&-8&0&0&0&0&24&-1&17&4&-6&20\cr
6&6&0&-3&-1&0&9972&0&6&2&0&0&0&0&0&18&-2&15&4&-7&21\cr
4&-8&0&-1&1&-2&-3018&0&-16&-2&-8&0&0&0&0&25&-1&20&3&-5&23\cr
5&2&0&-2&0&-2&852&0&10&-6&16&0&0&0&0&26&0&17&4&-7&24\cr
-2&-12&8&-10&2&-8&8&-4&-51&26&-76&0&13&0&10&4&4&-7&-2&7&-20\cr
0&-24&6&-9&2&-8&24&-12&-45&17&-56&0&15&3&6&0&1&-4&-2&5&-12\cr
0&-9&3&-3&1&-2&24&-12&-30&23&-54&0&6&1&2&6&-4&8&0&1&-2\cr
0&-12&6&-9&3&-12&36&-18&-54&23&-64&0&18&-1&14&18&-3&15&-2&3&2\cr
-4&-15&1&-8&-1&-3&4&-2&-51&25&-71&0&17&-1&13&2&-5&9&8&-14&26\cr
2&-12&4&-14&5&-16&-8&4&-39&22&-61&0&11&-4&15&8&-2&10&0&0&8\cr
0&-3&3&-3&1&-4&48&-24&-39&11&-43&0&21&1&15&24&-7&28&2&-7&30\cr
0&-15&3&-12&4&-15&0&0&-48&23&-68&0&18&1&10&6&-1&9&-4&5&2\cr}
$$
Rather than two, we use the first $10$ Fourier coefficients and  hence obtain a grossly overdetermined 
linear system of $80$ equations in $21$ unknowns. 
As expected, the solution space has dimension $6$. 
In this way we obtain six independent quadrics $\sum_{1\le i\le j\le 6}a_{ij}x_ix_j$ with
coefficients in~$\QQ$. By means of a linear change of variables and by
replacing the quadrics by suitable linear combinations, we obtain equations that have very small 
coefficients and have good reduction modulo primes different from~$17$. Here we use the LLL-algorithm as in~[\MER].
The independent  quadrics $q_1,\ldots, q_6$ we obtained, are listed below. They cut out a genus $6$ curve, which must be~$X_{\rm ns}^+(17)$.
$$
\eqalign{q_1=&-3x_1x_2+x_1x_3+x_1x_4+x_1x_5+x_2x_3+2x_2x_4+x_2x_5-x_2x_6-2x_3^2+\cr
&+2x_3x_4+2x_3x_5+x_3x_6+x_4x_5-x_4x_6+x_5^2-x_5x_6,\cr
q_2=&\ \ x_1x_2-2x_1x_3-2x_1x_4+x_1x_6+x_2x_5+2x_2x_6-x_3x_4-2x_3x_5+x_4^2-\cr
&-x_4x_5+x_4x_6-2x_5^2+x_6^2,\cr
q_3=&\ \ 3x_1^2+3x_1x_2+x_1x_3-x_1x_4+x_1x_6+x_2x_3-x_2x_4+x_2x_5+2x_2x_6+ x_3^2-\cr
&-x_3x_4-x_4^2-x_4x_5-x_4x_6+x_5^2+2x_5x_6,\cr
q_4=&\ \ 2x_1^2+2x_1x_2-2x_1x_3+x_1x_4-2x_1x_5+x_1x_6-x_2x_3-x_2x_5+3x_2x_6-x_3^2+\cr
&+3x_3x_4-3x_3x_5-x_4^2-x_4x_5+2x_5^2-x_5x_6+x_6^2,\cr
q_5=&\ \ x_1x_2+5x_1x_3+2x_1x_4-x_1x_5+x_2^2+3x_2x_3+2x_2x_4-x_2x_5-x_3^2+2x_3x_4-\cr
& -3x_3x_5+x_4^2+3x_4x_6-x_5^2-2x_5x_6-x_6^2,\cr
q_6=&-3x_1x_2+x_1x_3-2x_1x_4+4x_1x_5-3x_1x_6-3x_2^2-2x_2x_3-5x_2x_4 +x_2x_5-\cr
&-x_2x_6+x_3^2+x_3x_4-3x_3x_5+x_4^2-2x_4x_5-2x_4x_6+x_5^2+3x_5x_6 -x_6^2.\cr}
$$
CM-points  or Heegner points are points on modular curves parametrizing elliptic curves  with
complex multiplication by  imaginary quadratic orders ${\cal O}\subset \CC$.
Only if ${\cal O}$ is one of the thirteen
quadratic orders of class number~$1$, the CM-points may give rise to rational points.  Since the prime $17$ is inert in the orders ${\cal O}$ of discriminant $-3,-7,-11,-12$, $-27$, $-28$ and $-163$,
there is for each of these orders  ${\cal O}$, a unique rational CM-point on the curve~$X_{\rm ns}^+(17)$.  We have determined the projective coordinates of these CM-points by  evaluating the Fourier series of the modular forms $g_i$ numerically  in suitable $\tau\in\HH$ for which $17\tau\in R$.
\medskip

\hbox{\vbox{\hsize 2in\noindent {\bf Table.} CM-points  on $X_{\rm ns}^+(17)$.\medskip\vbox{\offinterlineskip
\hrule\halign{&\vrule#&\strut\quad\hfil#\quad\cr
height3pt
&\omit&&\omit&\cr
&discriminant&&CM-point&\cr
height3pt
&\omit&&\omit&\cr
\noalign{\hrule}
height2pt
&\omit&&\omit&\cr
&$-3$&&$(2 :-2 :-1 : 3 :-2 : 1)$&\cr
&$-7$&&$(-6 :-2 :-4 : 1 :-3 : 13)$&\cr
&$-11$&&$(3 : 1 : 2 :-9 :-7 : 2) $&\cr
&$-12$&&$(-4 : 10 : 3 :-5 :-2 : 3) $&\cr
&$-27$&&$ (2 :-5 :-10 :-6 : 1 : 7) $&\cr
&$-28$&&$(0 : 0 : 0 : 1 : 1 : 1) $&\cr
&$-163$&&$(-7 : 9 : 35 : 21 : 5 : 1)$&\cr
height3pt
&\omit&&\omit&\cr
}\hrule}}}
\medskip\noindent
A short computer calculation revealed that there are no rational points $(x_1:x_2:x_3:x_4:x_5:x_6)$ on $X_{\rm ns}^+(17)$ with $x_i\in\ZZ$ and $|x_i|<10\,000$,
other than the seven CM-points listed in the table.

\beginsection 7.  Level $19$ and $23$.

In this section we present quadrics that cut out the modular curves $X_{\rm ns}^+(19)$  and $X_{\rm ns}^+(23)$. They were obtained by the method explained in the previous section.

The modular curve $X_{\rm ns}^+(19)$ has genus $8$. Its canonical embedding in ${\bf P}_7$ is cut out by  fifteen quadrics. These are listed in Table 1.  Here the rows contain the coefficients of  the $36$ monomials $x_ix_j$ with $1\le i\le j\le 8$
in lexicographic order. Each column corresponds to the equation of a quadric in~$\PP_7$. 

\medskip\noindent{\bf Table 1.}
\medskip
{\eightpoint
$$
\matrix{
 -1 & 1 & 0 & 0 & -1 & 0 & 0 & -1 & 1 & 1 & 0 & 2 & 0 & 1 & 0 \cr
 1 & 0 & 0 & 0 & 1 & 0 & 0 & -1 & -2 & -1 & 1 & 0 & 3 & -2 & 4 \cr
 0 & -1 & -1 & 0 & 1 & 0 & 1 & 2 & -1 & 0 & 1 & -3 & 1 & 0 & 1 \cr
 0 & 1 & -1 & 0 & -1 & 1 & 0 & 1 & -2 & 0 & 1 & 2 & 1 & 0 & 0 \cr
 -1 & -1 & 1 & -2 & 1 & 0 & 0 & -1 & 1 & -1 & 0 & 1 & 1 & 2 & 2 \cr
 0 & 0 & 0 & -1 & 2 & 0 & -1 & 0 & -1 & -1 & 0 & -2 & 0 & 0 & 1 \cr
 0 & 1 & 0 & -1 & 0 & 0 & 0 & 0 & -1 & -1 & -1 & -1 & 0 & 0 & -1 \cr
 0 & 0 & 0 & 0 & 1 & -1 & 1 & -1 & 0 & 0 & 0 & 1 & 2 & -1 & -1 \cr
 0 & 0 & 1 & 1 & 0 & 0 & 0 & 0 & 0 & 0 & -1 & 0 & 0 & 0 & 1 \cr
 1 & -1 & 1 & 0 & -1 & 0 & -1 & 2 & 1 & -1 & 0 & 0 & -1 & 0 & 0 \cr
 0 & -1 & 1 & 0 & 1 & 0 & 0 & 0 & -1 & 1 & 0 & -1 & 2 & -1 & 0 \cr
 0 & 1 & -1 & 0 & 0 & 1 & 0 & -1 & 0 & 1 & 0 & 0 & 1 & -1 & -1 \cr
 0 & 0 & 0 & 0 & 0 & 2 & -1 & 0 & 1 & 2 & -1 & 0 & 2 & -1 & -1 \cr
 1 & 0 & -1 & -1 & -1 & 1 & 1 & 1 & -1 & 0 & 0 & 0 & 0 & 1 & -2 \cr
 0 & 0 & 1 & 1 & 0 & 0 & 0 & 0 & 0 & -1 & 1 & 1 & 0 & 1 & -1 \cr
 0 & -2 & 0 & 0 & 0 & 0 & 0 & 2 & 0 & 0 & -1 & -1 & -1 & 0 & 0 \cr
 0 & 1 & 1 & 0 & -1 & -1 & 0 & -2 & 2 & -1 & 0 & 0 & 0 & 1 & -1 \cr
 1 & -1 & 0 & -2 & 1 & -1 & 0 & 1 & 1 & 1 & 2 & -1 & 0 & 0 & -1 \cr
 1 & -1 & 1 & -1 & 0 & -1 & -1 & -2 & 2 & 0 & -1 & 0 & 0 & 0 & -1 \cr
 1 & -1 & 0 & -1 & -1 & 1 & 0 & 1 & -1 & 0 & -1 & 0 & 0 & 0 & 0 \cr
 1 & 0 & 1 & 0 & -2 & 2 & -2 & 1 & -1 & 0 & 0 & -1 & 0 & -1 & 0 \cr
 1 & 0 & 1 & 0 & 0 & 1 & 0 & 1 & -2 & -1 & 1 & 1 & 1 & -1 & 1 \cr
 0 & 0 & -1 & 0 & -1 & 0 & 0 & -1 & 0 & -1 & 0 & 2 & -1 & -1 & -3 \cr
 1 & -1 & 1 & 0 & 0 & 1 & 0 & 0 & 0 & 0 & -1 & 0 & 0 & -3 & 0 \cr
 -1 & 2 & 0 & 0 & -1 & 1 & 0 & 0 & 0 & 1 & -1 & -1 & -2 & 0 & -1 \cr
 1 & 0 & 1 & 0 & 1 & 0 & 0 & -1 & 0 & -1 & 0 & 1 & 1 & -1 & 1 \cr
 -1 & 0 & -1 & -1 & 0 & 0 & 0 & -1 & 1 & 0 & 0 & 0 & -1 & 0 & -2 \cr
 -1 & 1 & -1 & 0 & 0 & -1 & 2 & -2 & 0 & -1 & 0 & 0 & -2 & -1 & -1 \cr
 0 & -2 & 1 & -1 & 1 & -2 & -1 & 1 & 0 & 1 & -1 & -1 & -2 & -1 & -1 \cr
 -1 & 0 & 0 & 1 & 1 & -1 & -1 & 0 & -1 & 1 & 1 & -1 & 0 & -3 & -1 \cr
 0 & 0 & 0 & 0 & 0 & 0 & 1 & 0 & 0 & 0 & 0 & 0 & 0 & -1 & 0 \cr
 0 & 0 & 0 & -1 & 0 & 0 & 0 & 2 & 0 & 1 & 1 & 0 & -2 & 0 & 0 \cr
 0 & 0 & 0 & 0 & 1 & 0 & -1 & 0 & 0 & 0 & 0 & -1 & -3 & -1 & 1 \cr
 1 & 0 & 0 & -1 & -1 & 0 & 0 & 0 & 1 & 0 & 1 & 1 & 0 & 0 & 0 \cr
 0 & 0 & -1 & -1 & -2 & 1 & 1 & 1 & 1 & 0 & -1 & 0 & -1 & 1 & 1 \cr
 0 & 0 & 0 & 0 & 1 & 0 & -1 & 0 & 0 & 0 & 0 & -1 & 0 & 0 & 2 \cr
}
$$}

\bigskip\noindent
The  prime $19$ is inert in the imaginary quadratic orders ${\cal O}$ of discriminant $-4,-7,-11,-16$, $-28$, $-43$ and $-163$.
For each order ${\cal O}$ there is a rational CM-point   on $X_{\rm ns}^+(19)$,  corresponding to an elliptic curve with complex multiplication by~${\cal O}$.  As in the previous section, the CM-points have been computed numerically.  They are the only rational points $(x_1:x_2:x_3:x_4:x_5:x_6:x_7:x_8)$ with $x_i\in\ZZ$ satisfying $|x_i|\le 10\,000$.

\medskip
\hbox{\vbox{\hsize 2in\noindent {\bf Table 2.} CM-points on $X_{\rm ns}^+(19)$.\medskip\vbox{\offinterlineskip
\hrule\halign{&\vrule#&\strut\quad\hfil#\quad\cr
height3pt
&\omit&&\omit&\cr
&discriminant&&CM-point&\cr
height3pt
&\omit&&\omit&\cr
\noalign{\hrule}
height2pt
&\omit&&\omit&\cr
&$-4$&&$(0 : 0 :-1 : 1 : 0 :-1 : 1 : 0) $&\cr
&$-7$&&$(2 :7 : -12 : -4 :3 :3 :10 : -4)$&\cr
&$-11$&&$(3 : 1 : 1 :-6 :-5 :-5 :-4 : 13)  $&\cr
&$-16$&&$(-2 : 12 :7 :-15 : 16 :-3 : 9 : 4) $&\cr
&$-28$&&$ (0 :1 : 0 : 0 :1 : -1 : 0 : 0)  $&\cr
&$-43$&&$(-10 : 3 : 3 : 1 : 4 :-15 : 7 : 1) $&\cr
&$-163$&&$(2 : 0 : 0 :-3 :-1 : 0 : 0 : 3) $&\cr
height3pt
&\omit&&\omit&\cr
}\hrule}}}

\bigskip\noindent
The modular curve $X_{\rm ns}^+(23)$ has genus $13$. Its canonical embedding in ${\bf P}_{12}$ is cut out by  $55$ quadrics. These are listed in Table 3.  Here the rows contain the coefficients of  the $78$ monomials $x_ix_j$ with $1\le i\le j\le 13$
in lexicographic order. Each column corresponds to the equation of a quadric in~$\PP_{12}$. 

\hfuzz 2pt
\medskip\noindent{\bf Table 3.}
{\eightpoint
$$
\def\quad{\hskip0ex\relax}
\def\-{%
  \setbox0=\hbox{-}%
  \vcenter{%
    \hrule width\wd0 height \the\fontdimen8\textfont3%
  }%
}
\matrix{
 0 & 1 & \-1 & 1 & 3 & 0 & 0 & 2 & 1 & 0 & \-1 & 0 & 0 & 0 & 1 & \-2 & 1 & 2 & 2 & 1 & 1 & 0 & 1 & 0 & 0 & 0 & 1 & \-1 & 1 & 3 & \-2 & 1 & \-1 & \-1 & 0 & 0 & \-1 & 0 & 0 & 1 & \-1 & 1 & \-1 & \-1 & 0 & 2 & 2 & \-1 & \-1 & 2 & 1 & 1 & 1 & \-1 & 0 \cr
 \-2 & 1 & \-1 & 2 & 2 & 1 & \-1 & 0 & \-2 & \-1 & \-3 & 0 & 0 & \-2 & 1 & 2 & 3 & \-1 & \-1 & \-1 & 0 & 1 & 3 & \-1 & \-1 & 0 & \-2 & 0 & \-1 & \-1 & 1 & 0 & 1 & \-1 & 2 & 1 & \-3 & 1 & 0 & 0 & 1 & 0 & \-1 & \-1 & \-1 & \-1 & \-1 & 1 & 0 & 2 & \-1 & \-2 & \-1 & \-1 & 2 \cr
 0 & 0 & 2 & \-3 & \-2 & 2 & \-2 & 0 & \-2 & 1 & 0 & 0 & 0 & \-2 & \-3 & \-1 & 2 & \-2 & \-2 & 0 & \-1 & 0 & 1 & \-1 & \-1 & 0 & \-1 & \-1 & \-3 & \-2 & 3 & 1 & 1 & 1 & \-1 & 1 & 0 & 0 & 2 & \-2 & 0 & 0 & 0 & 1 & 1 & \-1 & \-1 & 0 & 0 & 0 & 0 & 0 & 0 & 0 & 0 \cr
 2 & 0 & 0 & 1 & \-1 & \-1 & 1 & 1 & 1 & 0 & 0 & 1 & 0 & 0 & 2 & 0 & \-1 & 1 & 1 & 1 & 1 & \-2 & \-1 & \-1 & 2 & \-1 & 1 & \-2 & 1 & 2 & \-2 & \-2 & \-1 & 0 & \-1 & \-3 & 3 & \-1 & \-1 & 2 & 0 & 0 & 0 & 2 & \-2 & 0 & 1 & 0 & 0 & 0 & 0 & 0 & 0 & 2 & 0 \cr
 \-1 & 0 & 0 & \-1 & 0 & 0 & \-1 & 0 & \-1 & 0 & 1 & \-1 & 0 & 0 & \-1 & \-1 & 1 & 0 & 0 & 0 & \-3 & 1 & \-1 & \-1 & \-2 & 1 & 0 & 1 & \-1 & \-1 & 1 & 1 & 0 & 1 & \-1 & 1 & \-1 & \-1 & 1 & \-1 & 0 & \-1 & 0 & \-1 & 1 & 0 & 0 & 1 & \-1 & 1 & 1 & 2 & 0 & \-1 & 0 \cr
 0 & 0 & 1 & 0 & 0 & \-2 & 2 & \-1 & \-3 & 0 & \-1 & \-1 & 1 & 2 & 3 & 2 & 1 & 1 & \-5 & 3 & \-3 & 2 & 1 & 0 & 0 & 2 & 1 & \-2 & 1 & \-1 & \-2 & 0 & \-2 & 1 & 0 & 1 & \-2 & 4 & 0 & 0 & \-3 & \-2 & 0 & 1 & 0 & \-1 & \-1 & 2 & 1 & \-1 & 1 & 0 & \-1 & \-1 & 0 \cr
 \-1 & \-2 & 1 & 0 & 3 & 2 & 0 & \-3 & 1 & 0 & \-2 & 0 & 0 & 0 & 1 & 1 & 2 & 0 & \-2 & 0 & 0 & 0 & 2 & \-1 & \-3 & 3 & 1 & \-2 & 2 & \-1 & \-1 & 1 & 1 & 1 & \-1 & 1 & \-2 & 1 & 1 & 3 & \-1 & 1 & 2 & 1 & 2 & \-1 & 0 & 1 & 2 & \-1 & 2 & \-1 & 1 & \-1 & 1 \cr
 0 & \-2 & 1 & 0 & \-1 & 2 & \-3 & 0 & 3 & 1 & \-2 & 2 & \-1 & \-5 & 1 & 3 & 1 & \-1 & 2 & \-2 & 3 & \-3 & 2 & 0 & 0 & \-2 & \-1 & 1 & 0 & 1 & 2 & 1 & 1 & \-1 & 0 & 0 & \-1 & \-2 & 1 & 2 & \-1 & 1 & 1 & 1 & 2 & 0 & 1 & 0 & 2 & 0 & 0 & \-2 & 0 & 2 & 1 \cr
 \-1 & 1 & 0 & 0 & 3 & \-3 & 0 & 3 & \-1 & 2 & \-1 & 0 & 1 & 1 & \-1 & 0 & 1 & 1 & 0 & 0 & 1 & 0 & 2 & \-3 & \-1 & \-1 & 0 & \-2 & 1 & 1 & \-1 & 5 & 0 & \-2 & 0 & 1 & 1 & 2 & 0 & 0 & \-2 & 2 & \-2 & \-2 & 2 & 2 & 0 & 2 & \-3 & \-1 & 0 & 1 & \-2 & 0 & \-1 \cr
 \-1 & \-1 & 1 & 0 & \-1 & 0 & \-4 & 0 & \-2 & \-1 & \-1 & \-2 & 1 & \-3 & 4 & 2 & 1 & \-4 & 0 & 1 & \-1 & \-1 & 0 & 1 & 0 & 4 & \-2 & \-2 & \-2 & \-4 & 0 & \-1 & 0 & 2 & 0 & \-1 & \-1 & \-2 & 0 & 0 & 0 & \-2 & \-2 & 4 & \-1 & \-3 & 0 & 3 & 1 & 0 & \-1 & \-1 & \-2 & 1 & 1 \cr
 0 & \-4 & \-1 & 2 & 0 & \-1 & 1 & 0 & 4 & \-3 & 1 & 2 & \-3 & \-1 & 2 & 0 & \-1 & 4 & 0 & \-1 & \-1 & 1 & 1 & \-1 & \-1 & \-2 & 1 & 2 & 1 & 2 & \-2 & \-2 & \-2 & 0 & \-1 & \-2 & \-1 & \-2 & \-1 & 3 & 0 & \-1 & \-1 & 0 & 2 & 3 & 1 & 1 & \-1 & 0 & 3 & \-2 & 1 & 1 & 1 \cr
 1 & \-1 & \-1 & 2 & \-2 & 1 & \-1 & 1 & 0 & 1 & 0 & \-1 & 1 & \-2 & 0 & \-1 & \-2 & \-1 & 0 & \-3 & 2 & \-2 & 1 & \-3 & 1 & \-1 & \-1 & 3 & \-2 & \-2 & \-2 & \-1 & 0 & 1 & 0 & \-2 & 3 & \-3 & \-1 & \-1 & 1 & \-1 & 1 & 2 & 2 & \-2 & 0 & \-1 & 0 & 1 & 0 & \-1 & 0 & 1 & \-2 \cr
 \-3 & \-2 & 0 & \-1 & 0 & 1 & \-1 & \-3 & 2 & \-1 & 1 & 2 & 0 & 1 & 1 & 4 & 2 & 0 & \-1 & \-2 & \-2 & 1 & 0 & 1 & 0 & \-1 & \-1 & 4 & 1 & \-2 & 2 & 1 & 1 & \-1 & 1 & 3 & \-3 & 0 & 3 & 1 & 0 & \-1 & 1 & \-4 & 1 & \-1 & \-1 & 1 & 2 & \-1 & 1 & 0 & 0 & 1 & 3 \cr
 0 & 2 & \-1 & 3 & \-1 & 0 & 0 & 2 & \-2 & \-1 & \-2 & \-1 & \-1 & \-3 & 1 & 1 & 0 & \-2 & 1 & 0 & 0 & \-1 & 1 & \-1 & 3 & \-1 & \-1 & \-1 & \-1 & 1 & 1 & \-2 & 0 & \-1 & 2 & \-2 & 1 & 0 & \-1 & \-1 & 2 & 0 & \-1 & 1 & \-3 & \-1 & \-1 & 1 & \-1 & 2 & \-4 & \-1 & \-1 & 0 & 1 \cr
 0 & 0 & 1 & 1 & \-1 & 1 & \-1 & 0 & \-1 & \-2 & 0 & 0 & 0 & \-1 & 1 & 2 & 1 & \-1 & \-1 & 0 & 0 & 0 & 1 & 2 & 2 & 2 & \-1 & \-2 & 0 & \-2 & 2 & \-1 & 1 & 1 & \-1 & \-2 & \-1 & 0 & 1 & \-1 & 0 & 0 & 0 & 0 & \-2 & \-1 & 0 & 1 & \-1 & 0 & \-1 & \-2 & \-1 & \-1 & 0 \cr
 \-2 & 0 & \-1 & \-1 & 2 & 0 & 0 & \-1 & 0 & 1 & 0 & \-1 & 0 & 1 & \-2 & \-1 & 0 & 0 & 0 & \-1 & 0 & 2 & 1 & \-2 & \-2 & 0 & 0 & 1 & 0 & 0 & 1 & 3 & 2 & 0 & 1 & 3 & 0 & 0 & 0 & 1 & 0 & 0 & 0 & \-1 & 3 & 0 & 0 & 0 & 2 & 0 & 1 & 2 & 1 & \-1 & 0 \cr
 2 & 0 & 1 & 1 & \-1 & 1 & 1 & 0 & 1 & 0 & \-2 & 2 & 0 & \-1 & 1 & 2 & \-1 & \-1 & \-1 & 0 & 2 & \-2 & \-1 & 1 & 1 & 1 & \-1 & \-1 & 1 & 0 & \-1 & \-2 & \-1 & 1 & \-1 & \-2 & \-1 & 0 & \-1 & 1 & 0 & 0 & 1 & 1 & \-1 & \-1 & 0 & 1 & \-1 & \-1 & 0 & \-2 & \-1 & 1 & 0 \cr
 0 & 0 & \-2 & 2 & 0 & 0 & 1 & \-1 & \-2 & 1 & \-3 & 1 & 2 & 2 & \-2 & 2 & 1 & 1 & \-3 & 0 & 1 & 1 & 1 & \-3 & \-1 & \-1 & 0 & 0 & 0 & 0 & \-2 & 2 & 0 & 0 & 2 & \-1 & 0 & 0 & \-1 & 2 & 3 & 1 & \-3 & 1 & 2 & \-2 & \-3 & 1 & 0 & 2 & 0 & \-1 & \-4 & \-2 & 0 \cr
 \-3 & 1 & \-2 & 2 & 2 & 1 & \-1 & \-1 & \-3 & 0 & \-3 & \-3 & \-1 & \-2 & \-1 & \-1 & 2 & \-1 & 0 & \-1 & \-2 & 1 & 2 & \-3 & 0 & 1 & 0 & \-1 & \-1 & \-1 & 1 & 1 & 1 & \-1 & 3 & 1 & \-1 & 0 & 1 & 0 & 2 & 0 & 0 & 0 & 0 & \-2 & \-2 & 1 & 2 & 3 & \-2 & 1 & 0 & \-1 & 2 \cr
 \-1 & 1 & \-1 & 1 & 3 & 2 & \-3 & 0 & \-1 & 1 & \-1 & \-2 & 1 & 0 & 0 & \-1 & 3 & 2 & 2 & 0 & 0 & \-1 & 0 & \-1 & \-1 & 1 & 1 & \-1 & \-1 & 0 & 0 & 2 & 0 & 0 & \-1 & 1 & \-2 & \-1 & 2 & 0 & \-1 & 2 & 1 & \-2 & \-1 & 0 & 2 & \-1 & 0 & 4 & 0 & 1 & 2 & \-1 & 0 \cr
 1 & 0 & 1 & \-1 & \-2 & \-2 & 2 & \-1 & 0 & 1 & 0 & 3 & 1 & 0 & 0 & 0 & 0 & 1 & \-1 & \-1 & \-1 & 0 & \-1 & 1 & 0 & \-1 & \-2 & 0 & 2 & \-1 & \-1 & \-1 & 0 & \-2 & 1 & \-2 & 0 & 0 & 0 & \-1 & \-2 & \-2 & 0 & 1 & \-2 & 0 & \-1 & \-1 & 1 & \-2 & 2 & 0 & \-2 & 2 & 2 \cr
 0 & 0 & 0 & 0 & 2 & 2 & \-2 & 1 & 1 & 0 & 0 & \-1 & \-1 & \-2 & 2 & 2 & 0 & \-2 & 0 & 0 & 2 & \-2 & 0 & \-2 & \-1 & 1 & 1 & 0 & \-2 & 0 & 1 & \-1 & 1 & 3 & \-1 & 3 & 1 & 0 & \-2 & 0 & 0 & 1 & 1 & \-1 & 0 & \-1 & 1 & 1 & \-1 & \-1 & \-1 & \-1 & 3 & 0 & \-1 \cr
 1 & 0 & 0 & \-1 & 0 & \-1 & 2 & 0 & 3 & \-1 & \-1 & 0 & \-4 & 0 & 0 & 0 & 0 & \-2 & \-1 & \-1 & \-1 & 1 & \-2 & \-2 & 1 & \-2 & 0 & \-1 & 2 & 4 & 0 & 0 & \-4 & 0 & 1 & 0 & \-2 & 2 & 0 & 4 & \-2 & 2 & 1 & \-1 & 0 & 0 & \-2 & 3 & \-1 & \-1 & 0 & \-3 & 1 & 3 & 2 \cr
 0 & \-2 & 2 & \-2 & \-1 & 1 & \-2 & 0 & 1 & \-2 & 1 & \-1 & \-4 & \-1 & 2 & 2 & \-1 & \-3 & 1 & \-2 & 0 & 1 & 1 & \-1 & 3 & 1 & \-2 & 0 & \-1 & \-1 & 1 & 0 & \-1 & 5 & \-3 & 0 & 0 & 2 & 0 & 0 & 0 & \-1 & 0 & 0 & 3 & 1 & 0 & 2 & \-1 & \-1 & \-1 & \-2 & 1 & 1 & \-1 \cr
 0 & \-1 & 1 & \-1 & \-1 & 0 & 0 & 1 & 1 & 0 & 1 & 0 & \-1 & \-1 & 1 & 0 & \-1 & \-1 & 1 & 2 & \-2 & \-2 & 0 & 0 & 0 & 0 & 1 & 1 & \-3 & 1 & 0 & 0 & \-1 & 0 & \-1 & 1 & 1 & 1 & 0 & \-1 & 0 & 0 & 0 & 1 & 2 & 1 & 1 & 1 & \-1 & \-1 & \-2 & 2 & 2 & 1 & 0 \cr
 0 & 0 & 0 & 0 & 0 & 0 & 0 & 0 & 0 & 0 & 0 & 0 & 0 & 0 & 0 & 0 & 0 & 0 & 0 & 0 & 0 & 0 & 0 & 1 & 1 & 1 & 0 & 0 & 0 & 0 & 0 & 0 & 0 & 0 & 0 & 0 & 0 & 0 & 0 & 0 & \-1 & 1 & 0 & 0 & 0 & 0 & 0 & 0 & 0 & 0 & 0 & 0 & 0 & 0 & 0 \cr
 0 & 0 & 0 & 0 & 0 & 0 & 0 & 0 & 0 & 0 & 0 & 0 & 0 & 0 & 0 & 0 & 0 & 0 & 0 & 0 & 1 & 0 & 1 & 1 & \-1 & \-1 & 0 & 0 & 0 & 0 & 0 & 0 & 0 & \-1 & 1 & 0 & 0 & 1 & 0 & 0 & 0 & 0 & 1 & 0 & 1 & 0 & 1 & 0 & \-1 & 0 & 0 & 0 & 0 & 0 & 0 \cr
 0 & 0 & 0 & 0 & 0 & 0 & 0 & 0 & 0 & 0 & 0 & 0 & 0 & 0 & 0 & 0 & 0 & 0 & 0 & 0 & 0 & 0 & 0 & 0 & 0 & 1 & 0 & 0 & \-1 & 0 & 0 & 0 & 0 & 0 & 0 & 0 & 0 & 0 & 0 & 0 & 0 & 0 & 0 & 0 & 0 & 1 & 0 & 0 & 0 & 0 & 0 & \-1 & 0 & 0 & 0 \cr
 0 & 0 & 0 & 0 & 0 & 0 & 0 & 0 & 0 & 0 & 0 & 0 & 0 & 0 & 0 & 0 & 0 & 0 & 0 & 0 & 0 & 0 & 0 & 0 & 0 & 0 & 1 & \-1 & 1 & 0 & \-1 & 0 & 0 & 0 & 1 & 0 & 0 & 0 & 0 & 0 & 1 & 0 & 0 & 0 & 0 & 0 & 0 & 0 & \-1 & 0 & 0 & 0 & \-1 & \-1 & 0 \cr
 0 & 0 & 0 & 0 & 0 & 0 & 0 & 0 & 0 & 0 & 0 & 0 & 0 & 0 & 0 & 0 & 0 & 0 & 0 & 0 & 0 & 0 & 0 & 0 & 0 & 0 & 1 & 0 & 0 & 0 & 0 & 0 & 0 & 1 & 0 & 0 & 0 & \-1 & 0 & 0 & 0 & 0 & 0 & 0 & 0 & 0 & 0 & 1 & 0 & 0 & 0 & 0 & 0 & 0 & 0 \cr
 0 & 0 & 0 & 0 & 0 & 0 & 0 & 0 & 0 & 0 & 0 & 0 & 0 & 0 & 0 & 0 & 0 & 0 & 0 & 0 & 0 & 0 & 0 & 0 & 0 & 0 & 0 & 1 & 0 & 0 & 0 & 0 & 0 & 0 & 0 & 0 & 0 & \-1 & 0 & 0 & 0 & 0 & 0 & 0 & 0 & 0 & 0 & 0 & 1 & 0 & 0 & 0 & 0 & 0 & 0 \cr
 0 & 0 & 0 & 0 & 0 & 0 & 0 & 0 & 0 & 0 & 0 & 0 & 0 & 0 & 0 & 0 & 0 & 0 & 0 & 0 & 0 & 0 & 0 & 0 & 0 & 0 & 0 & 0 & 1 & 1 & \-1 & \-1 & 0 & 0 & 1 & 0 & 0 & 0 & 0 & 0 & 0 & 0 & 0 & 1 & 0 & 0 & 1 & 0 & 0 & 1 & 0 & 0 & 0 & 0 & 0 \cr}
 $$
 
 $$\def\quad{\hskip0.01ex\relax}
 \def\quad{\hskip0ex\relax}
\def\-{%
  \setbox0=\hbox{-}%
  \vcenter{%
    \hrule width\wd0 height \the\fontdimen8\textfont3%
  }%
}
 \matrix{
 0 & 0 & 0 & 0 & 0 & 0 & 0 & 0 & 0 & 0 & 0 & 0 & 0 & 0 & 0 & 0 & 0 & 0 & 0 & 0 & 0 & 0 & 0 & 0 & 0 & 0 & 0 & 0 & 0 & 1 & 1 & 0 & 1 & \-1 & 0 & 0 & 1 & 0 & 0 & 0 & 0 & 0 & 0 & 0 & 0 & 0 & 0 & 0 & 0 & 0 & 1 & 0 & 0 & 0 & 0 \cr
 0 & 0 & 0 & 0 & 0 & 0 & 0 & 0 & 0 & 0 & 0 & 0 & 0 & 0 & 0 & 0 & 0 & 0 & 0 & 0 & 0 & 0 & 0 & 0 & 0 & 0 & 0 & 0 & 0 & 0 & 0 & 1 & 0 & 0 & 0 & 0 & \-1 & 1 & 0 & 0 & 0 & 0 & 0 & 0 & 0 & 0 & 0 & 0 & 0 & 0 & 0 & 0 & 0 & 0 & 0 \cr
 0 & 0 & 0 & 0 & 0 & 0 & 0 & 0 & 0 & 0 & 0 & 0 & 0 & 0 & 0 & 0 & 0 & 0 & 0 & 0 & 0 & 0 & 0 & 0 & 0 & 0 & 0 & 0 & 0 & 0 & 0 & 0 & 1 & 1 & \-1 & 0 & 0 & 0 & 0 & 0 & 0 & 0 & 0 & 0 & 0 & 0 & 0 & 0 & 0 & 0 & 0 & 0 & 0 & 0 & 0 \cr
 0 & 0 & 0 & 0 & 0 & 0 & 0 & 0 & 0 & 0 & 0 & 0 & 0 & 0 & 0 & 0 & 0 & 0 & 0 & 0 & 0 & 0 & 0 & 0 & 0 & 0 & 0 & 0 & 0 & 0 & 0 & 0 & 0 & 0 & 1 & 0 & 0 & 0 & 0 & 0 & 0 & 0 & 0 & 0 & 0 & 0 & 0 & 0 & 0 & 0 & 0 & 0 & 0 & 0 & 0 \cr
 0 & 0 & 0 & 0 & 0 & 0 & 0 & 0 & 0 & 0 & 0 & 0 & 0 & 0 & 0 & 0 & 0 & 0 & 0 & 0 & 0 & 0 & 0 & 0 & 0 & 0 & 0 & 0 & 0 & 0 & 0 & 0 & 0 & 0 & 0 & 0 & 0 & 0 & 1 & 0 & 0 & 0 & 0 & 0 & 0 & 0 & 0 & 0 & 0 & 0 & 0 & 0 & 0 & 0 & 0 \cr
 0 & 0 & 0 & 0 & 0 & 0 & 0 & 0 & 0 & 0 & 0 & 0 & 0 & 0 & 0 & 0 & 0 & 0 & 0 & 0 & 0 & 0 & 0 & 0 & 0 & 0 & 0 & 0 & 0 & 0 & 0 & 0 & 0 & 0 & 0 & 0 & 1 & 0 & 0 & 0 & 0 & 0 & 0 & 0 & 0 & 0 & 0 & 0 & 0 & 0 & 0 & 0 & 0 & 0 & 0 \cr
 0 & 0 & 0 & 0 & 0 & 0 & 0 & 0 & 0 & 0 & 0 & 0 & 0 & 0 & 0 & 0 & 0 & 0 & 0 & 0 & 0 & 0 & 0 & 0 & 0 & 0 & 0 & 0 & 0 & 0 & 0 & 0 & 0 & 0 & 0 & 1 & 0 & 0 & 1 & 0 & 0 & 0 & 0 & 0 & 0 & 0 & 0 & 0 & 0 & 0 & 0 & 0 & 0 & 0 & 0 \cr
 0 & 0 & 0 & 0 & 0 & 0 & 0 & 0 & 0 & 0 & 0 & 0 & 0 & 0 & 0 & 0 & 0 & 0 & 0 & 0 & 0 & 0 & 0 & 0 & 0 & 0 & 0 & 0 & 0 & 0 & 0 & 0 & 0 & 0 & 0 & 0 & 0 & 0 & 0 & 1 & 0 & 0 & 0 & 0 & 0 & 0 & 0 & 0 & 0 & 0 & 0 & 0 & 0 & 0 & 0 \cr
 0 & 0 & 0 & 0 & 0 & 0 & 0 & 0 & 0 & 0 & 0 & 0 & 0 & 0 & 0 & 0 & 0 & 0 & 0 & 0 & 0 & 0 & 0 & 0 & 0 & 0 & 0 & 0 & 0 & 0 & 0 & 0 & 0 & 0 & 0 & 0 & 0 & 0 & 0 & 0 & 1 & 0 & 0 & 0 & 0 & 0 & 0 & 0 & 0 & 0 & 0 & 0 & 0 & 0 & 0 \cr
 0 & 0 & 0 & 0 & 0 & 0 & 0 & 0 & 0 & 0 & 0 & 0 & 0 & 0 & 0 & 0 & 0 & 0 & 0 & 0 & 0 & 0 & 0 & 0 & 0 & 0 & 0 & 0 & 0 & 0 & 0 & 0 & 0 & 0 & 0 & 0 & 0 & 0 & 0 & 0 & 0 & 1 & 0 & 0 & 0 & 0 & 0 & 0 & 0 & 0 & 0 & 0 & 0 & 0 & 0 \cr
 0 & 0 & 0 & 0 & 0 & 0 & 0 & 0 & 0 & 0 & 0 & 0 & 0 & 0 & 0 & 0 & 0 & 0 & 0 & 0 & 0 & 0 & 0 & 0 & 0 & 0 & 0 & 0 & 0 & 0 & 0 & 0 & 0 & 0 & 0 & 0 & 0 & 0 & 0 & 0 & 0 & 0 & 1 & 0 & 0 & 0 & 0 & 0 & 0 & 0 & 0 & 0 & 0 & 0 & 0 \cr
 0 & 0 & 0 & 0 & 0 & 0 & 0 & 0 & 0 & 0 & 0 & 0 & 0 & 0 & 0 & 0 & 0 & 0 & 0 & 0 & 0 & 0 & 0 & 0 & 0 & 0 & 0 & 0 & 0 & 0 & 0 & 0 & 0 & 0 & 0 & 0 & 0 & 0 & 0 & 0 & 0 & 0 & 0 & 0 & 1 & 0 & 0 & 0 & 0 & 0 & 0 & 0 & 0 & 0 & 0 \cr
 0 & 0 & 0 & 0 & 0 & 0 & 0 & 0 & 0 & 0 & 0 & 0 & 0 & 0 & 0 & 0 & 0 & 0 & 0 & 0 & 0 & 0 & 0 & 0 & 0 & 0 & 0 & 0 & 0 & 0 & 0 & 0 & 0 & 0 & 0 & 0 & 0 & 0 & 0 & 0 & 0 & 0 & 0 & 1 & 0 & 0 & 0 & 0 & 0 & 0 & 0 & 0 & 0 & 0 & 0 \cr
 \-4 & 0 & \-2 & 2 & 4 & 0 & 0 & \-2 & 0 & \-2 & 0 & \-2 & 0 & 2 & 0 & 2 & 0 & 0 & 0 & \-2 & \-1 & 4 & 1 & \-1 & \-2 & 2 & 0 & 2 & 3 & 0 & 2 & 2 & 4 & 1 & 1 & 3 & \-3 & \-1 & 0 & 1 & 2 & \-1 & 0 & \-3 & 1 & 0 & 0 & 2 & 2 & 1 & 1 & 0 & 0 & \-4 & 0 \cr
 0 & 0 & 0 & 0 & 0 & 0 & 0 & 0 & 0 & 0 & 0 & 0 & 0 & 0 & 0 & 0 & 0 & 0 & 0 & 0 & 0 & 0 & 0 & 0 & 0 & 0 & 0 & 0 & 0 & 0 & 0 & 0 & 0 & 0 & 0 & 0 & 0 & 0 & 0 & 0 & 0 & 0 & 0 & 0 & 0 & 1 & 0 & 0 & 0 & 0 & 0 & 0 & 0 & 0 & 0 \cr
 0 & 0 & 0 & 0 & 0 & 0 & 0 & 0 & 0 & 0 & 0 & 0 & 0 & 0 & 0 & 0 & 0 & 0 & 0 & 0 & 0 & 0 & 0 & 0 & 0 & 0 & 0 & 0 & 0 & 0 & 0 & 0 & 0 & 0 & 0 & 0 & 0 & 0 & 0 & 0 & 0 & 0 & 0 & 0 & 0 & 0 & 1 & 0 & 0 & 0 & 0 & 0 & 0 & 0 & 0 \cr
 0 & 0 & 0 & 0 & 0 & 0 & 0 & 0 & 0 & 0 & 0 & 0 & 0 & 0 & 0 & 0 & 0 & 0 & 0 & 0 & 0 & 0 & 0 & 0 & 0 & 0 & 0 & 0 & 0 & 0 & 0 & 0 & 0 & 0 & 0 & 0 & 0 & 0 & 0 & 0 & 0 & 0 & 0 & 0 & 0 & 0 & 0 & 1 & 0 & 0 & 0 & 0 & 0 & 0 & 0 \cr
 0 & 0 & 0 & 0 & 0 & 0 & 0 & 0 & 0 & 0 & 0 & 0 & 0 & 0 & 0 & 0 & 0 & 0 & 0 & 0 & 0 & 0 & 0 & 0 & 0 & 0 & 0 & 0 & 0 & 0 & 0 & 0 & 0 & 0 & 0 & 0 & 0 & 0 & 0 & 0 & 0 & 0 & 0 & 0 & 0 & 0 & 0 & 0 & 1 & 0 & 0 & 0 & 0 & 0 & 0 \cr
 0 & 0 & 0 & 0 & 0 & 0 & 0 & 0 & 0 & 0 & 0 & 0 & 0 & 0 & 0 & 0 & 0 & 0 & 0 & 0 & 0 & 0 & 0 & 0 & 0 & 0 & 0 & 0 & 0 & 0 & 0 & 0 & 0 & 0 & 0 & 0 & 0 & 0 & 0 & 0 & 0 & 0 & 0 & 0 & 0 & 0 & 0 & 0 & 0 & 1 & 0 & 0 & 0 & 0 & 0 \cr
 0 & 0 & 0 & 0 & 0 & 0 & 0 & 0 & 0 & 0 & 0 & 0 & 0 & 0 & 0 & 0 & 0 & 0 & 0 & 0 & 0 & 0 & 0 & 0 & 0 & 0 & 0 & 0 & 0 & 0 & 0 & 0 & 0 & 0 & 0 & 0 & 0 & 0 & 0 & 0 & 0 & 0 & 0 & 0 & 0 & 0 & 0 & 0 & 0 & 0 & 1 & 0 & 0 & 0 & 0 \cr
 0 & 0 & 0 & 0 & 0 & 0 & 0 & 0 & 0 & 0 & 0 & 0 & 0 & 0 & 0 & 0 & 0 & 0 & 0 & 0 & 0 & 0 & 0 & 0 & 0 & 0 & 0 & 0 & 0 & 0 & 0 & 0 & 0 & 0 & 0 & 0 & 0 & 0 & 0 & 0 & 0 & 0 & 0 & 0 & 0 & 0 & 0 & 0 & 0 & 0 & 0 & 1 & 0 & 0 & 0 \cr
 0 & 0 & 0 & 0 & 0 & 0 & 0 & 0 & 0 & 0 & 0 & 0 & 0 & 0 & 0 & 0 & 0 & 0 & 0 & 0 & 0 & 0 & 0 & 0 & 0 & 0 & 0 & 0 & 0 & 0 & 0 & 0 & 0 & 0 & 0 & 0 & 0 & 0 & 0 & 0 & 0 & 0 & 0 & 0 & 0 & 0 & 0 & 0 & 0 & 0 & 0 & 0 & 1 & 0 & 0 \cr
 0 & 0 & 0 & 0 & 0 & 0 & 0 & 0 & 0 & 0 & 0 & 0 & 0 & 0 & 0 & 0 & 0 & 0 & 0 & 0 & 0 & 0 & 0 & 0 & 0 & 0 & 0 & 0 & 0 & 0 & 0 & 0 & 0 & 0 & 0 & 0 & 0 & 0 & 0 & 0 & 0 & 0 & 0 & 0 & 0 & 0 & 0 & 0 & 0 & 0 & 0 & 0 & 0 & 1 & 0 \cr
 0 & 0 & 0 & 0 & 0 & 0 & 0 & 0 & 0 & 0 & 0 & 0 & 0 & 0 & 0 & 0 & 0 & 0 & 0 & 0 & 0 & 0 & 0 & 0 & 0 & 0 & 0 & 0 & 0 & 0 & 0 & 0 & 0 & 0 & 0 & 0 & 0 & 0 & 0 & 0 & 0 & 0 & 0 & 0 & 0 & 0 & 0 & 0 & 0 & 0 & 0 & 0 & 0 & 0 & 1 \cr
 0 & 1 & \-2 & 0 & 3 & 1 & \-1 & \-1 & 0 & \-1 & 1 & 0 & 0 & 2 & 0 & 2 & 3 & 0 & 1 & 0 & 0 & 0 & \-1 & \-1 & \-2 & 2 & 1 & 0 & \-2 & \-1 & 1 & 2 & 1 & 2 & \-1 & 1 & 0 & 2 & \-2 & 1 & 2 & 4 & 0 & \-4 & \-1 & 0 & \-2 & \-1 & \-2 & 0 & \-1 & \-1 & 0 & \-2 & 0 \cr
 \-2 & 3 & \-3 & 3 & 4 & 1 & \-1 & 1 & \-1 & 2 & \-1 & \-4 & 2 & 3 & \-3 & 0 & 2 & 1 & 1 & \-2 & 1 & 0 & \-2 & \-3 & \-1 & \-1 & 2 & 0 & 3 & 2 & 0 & 5 & 1 & \-1 & 1 & 3 & \-2 & 0 & 2 & 1 & 0 & 5 & 1 & \-4 & 0 & \-2 & 0 & 1 & 0 & 5 & \-1 & 1 & 1 & \-4 & \-2 \cr
 0 & 1 & 3 & 1 & \-1 & 2 & 1 & 0 & \-4 & 2 & \-1 & \-2 & \-2 & 0 & \-3 & 1 & 0 & 1 & \-1 & 1 & \-1 & 1 & \-3 & \-1 & 1 & 1 & 5 & \-4 & 3 & 0 & 2 & 0 & \-2 & 3 & \-1 & 3 & \-2 & 1 & \-1 & 0 & 0 & 3 & 3 & 0 & \-2 & 0 & 0 & 2 & 0 & 0 & \-3 & \-2 & 2 & \-3 & 3 \cr
 0 & 1 & \-4 & 0 & 4 & \-3 & 1 & \-1 & 4 & \-2 & 3 & \-2 & 1 & 2 & \-1 & \-4 & \-1 & \-1 & 1 & 0 & \-1 & 3 & \-1 & \-1 & \-4 & 1 & \-3 & 2 & 0 & 1 & 1 & 3 & 2 & \-2 & 1 & \-1 & 1 & \-3 & \-1 & 1 & 0 & 2 & \-1 & \-3 & \-1 & 2 & 0 & \-1 & \-1 & \-1 & 2 & 1 & 0 & 1 & \-2 \cr
 3 & 0 & 1 & \-2 & 1 & 0 & 1 & \-2 & 3 & \-1 & 1 & 0 & \-2 & 1 & 2 & 2 & 1 & 1 & \-2 & 3 & 0 & 1 & \-1 & 0 & \-3 & 0 & 1 & \-2 & 0 & 4 & 2 & 0 & \-1 & 2 & \-2 & 1 & \-2 & 2 & 0 & 2 & \-3 & 2 & 3 & \-1 & \-1 & 2 & 2 & 0 & 0 & \-1 & 2 & \-3 & 3 & 1 & 0 \cr
 \-1 & \-2 & 0 & \-2 & 0 & 1 & \-2 & \-4 & \-2 & \-2 & 0 & 1 & \-3 & 2 & \-4 & 1 & 1 & 3 & \-1 & 0 & \-2 & 3 & 0 & 0 & \-1 & 1 & 0 & \-2 & \-2 & \-2 & 3 & 1 & 0 & 1 & \-1 & \-1 & \-1 & 0 & 1 & 0 & 4 & 1 & \-1 & \-2 & 2 & 3 & \-2 & 1 & \-3 & 1 & 1 & 0 & \-2 & 0 & 0 \cr
 \-1 & 0 & 0 & 2 & 0 & 0 & 0 & 2 & 0 & \-2 & \-1 & 0 & 2 & 3 & 2 & 2 & 1 & \-1 & 0 & 0 & 1 & \-1 & 2 & 1 & 2 & \-1 & \-1 & \-1 & 1 & 1 & \-2 & 1 & 0 & \-1 & \-2 & \-3 & \-1 & 2 & 3 & 1 & 1 & \-1 & \-3 & 1 & 2 & \-2 & 0 & 3 & \-1 & 3 & \-1 & 0 & \-2 & \-2 & \-2 \cr
 \-1 & \-1 & 0 & 0 & 1 & 1 & \-1 & \-2 & \-1 & 0 & \-1 & \-1 & 0 & 0 & 0 & \-1 & 2 & 1 & \-1 & 0 & \-2 & 0 & 0 & 0 & \-1 & 2 & 1 & \-1 & 0 & \-2 & \-1 & 0 & \-1 & 0 & 0 & 0 & \-2 & 0 & 1 & 1 & 0 & 0 & 1 & 0 & 0 & \-1 & \-1 & 0 & 1 & 1 & 1 & 0 & 0 & 0 & 1 \cr
 \-1 & \-1 & 2 & 1 & \-2 & 3 & \-2 & \-1 & \-3 & 2 & \-3 & 0 & 1 & \-4 & \-1 & 1 & 1 & 0 & \-1 & 0 & 0 & \-1 & 3 & 0 & 1 & 0 & 0 & 0 & 0 & \-2 & 1 & \-1 & 1 & 0 & 1 & 1 & \-1 & \-2 & 1 & \-1 & 1 & \-1 & 1 & 3 & 0 & \-1 & 0 & 0 & 2 & 2 & \-1 & 0 & \-1 & \-1 & 1 \cr
 \-1 & 0 & 1 & 3 & 2 & \-1 & 1 & 0 & \-1 & 2 & \-1 & \-1 & 2 & 1 & 2 & 3 & 0 & 1 & 1 & 1 & \-1 & \-1 & 0 & \-3 & \-1 & 2 & 3 & \-1 & 4 & 1 & \-1 & 2 & 0 & 1 & \-1 & 1 & 0 & 0 & \-1 & 2 & \-1 & 0 & 3 & 0 & 0 & 0 & 3 & 3 & 0 & \-1 & \-1 & 1 & 0 & \-1 & 0 \cr
 \-1 & \-1 & \-1 & \-2 & 1 & \-1 & 0 & \-1 & 2 & \-2 & 0 & \-1 & \-1 & 0 & \-1 & \-1 & \-2 & \-2 & \-4 & 0 & 0 & 2 & 1 & 0 & \-3 & 0 & \-3 & 0 & \-2 & 0 & 1 & 1 & 3 & \-1 & 1 & 1 & \-1 & \-1 & 1 & \-1 & 0 & \-1 & \-3 & 1 & 3 & 0 & \-2 & 2 & 1 & \-1 & 2 & 0 & 0 & 0 & \-1 \cr
 \-2 & \-4 & \-2 & 4 & 1 & 0 & 0 & \-4 & 0 & \-3 & 0 & 1 & \-1 & 0 & 1 & 1 & 0 & 5 & 1 & \-1 & \-3 & 2 & 1 & \-1 & \-3 & 2 & 1 & 1 & 1 & \-2 & \-1 & \-2 & 1 & 1 & 0 & \-3 & \-1 & \-3 & \-1 & 3 & 4 & \-2 & 0 & 1 & 0 & 2 & 1 & 1 & 0 & 2 & 2 & \-1 & \-1 & \-1 & 2 \cr
 \-1 & \-3 & 1 & \-1 & \-1 & \-1 & \-2 & \-2 & \-1 & \-2 & 0 & \-2 & \-1 & 2 & 0 & \-1 & 0 & 3 & \-4 & 1 & \-3 & 3 & 2 & \-1 & 1 & \-1 & \-1 & \-2 & 0 & \-1 & \-1 & 1 & \-2 & 0 & \-2 & \-2 & \-3 & \-1 & 4 & 0 & 0 & \-3 & \-1 & 1 & 3 & 2 & 0 & 3 & \-1 & 3 & 3 & 0 & \-2 & 2 & \-1 \cr
 1 & 1 & \-1 & 0 & \-1 & 2 & 1 & 0 & 1 & 1 & 0 & 2 & 0 & 0 & \-2 & \-1 & 0 & 2 & \-1 & \-1 & 1 & \-1 & \-1 & 0 & 1 & \-4 & 0 & 2 & 0 & 2 & 0 & \-2 & 0 & \-2 & 2 & 1 & 0 & 0 & 1 & \-2 & 0 & 1 & 2 & \-2 & \-1 & \-1 & 0 & \-3 & \-1 & 2 & 1 & 1 & 2 & 0 & 1 \cr
 1 & \-2 & 2 & \-1 & \-2 & 1 & 0 & \-2 & 0 & 0 & 0 & 2 & \-1 & \-2 & 0 & 2 & \-1 & 0 & \-1 & 1 & 0 & 0 & 1 & 0 & \-1 & 0 & 0 & 0 & \-1 & \-1 & 2 & \-2 & 1 & 2 & \-1 & 0 & 1 & \-1 & \-2 & 0 & 1 & \-1 & 1 & 2 & 0 & 1 & 0 & 0 & 0 & \-2 & 0 & \-1 & 0 & 1 & 1 \cr
 \-1 & \-1 & 0 & 2 & 0 & 0 & \-2 & 1 & 0 & 1 & \-1 & 0 & 3 & \-1 & 3 & 3 & 1 & 0 & 1 & 0 & 1 & \-2 & 3 & 0 & 1 & 0 & 0 & 2 & 1 & 0 & \-1 & 2 & 1 & \-1 & 0 & 0 & 0 & \-1 & 1 & 1 & \-1 & \-1 & 0 & 1 & 2 & \-1 & 2 & 1 & 1 & 1 & 0 & 1 & \-2 & \-1 & \-2 \cr
 2 & \-2 & 0 & \-2 & 1 & \-1 & 1 & \-1 & 3 & 0 & 0 & 3 & \-1 & 2 & 1 & 1 & \-2 & 2 & \-2 & 0 & 2 & 0 & \-1 & \-1 & \-2 & \-1 & \-1 & 0 & \-1 & 1 & \-3 & 0 & \-1 & 2 & \-1 & 0 & 0 & 1 & \-2 & 1 & \-1 & 0 & \-2 & 0 & 2 & 2 & 0 & 0 & \-1 & \-2 & 3 & \-1 & 1 & 2 & 0 \cr
 \-3 & \-3 & 0 & 2 & \-2 & 1 & \-3 & 0 & \-2 & \-1 & 0 & 0 & 0 & \-3 & \-1 & 0 & 1 & 2 & 1 & \-1 & \-2 & 0 & 2 & 0 & 0 & 0 & 1 & 1 & \-1 & \-3 & 1 & \-1 & 1 & \-1 & 0 & \-1 & 0 & \-4 & 1 & 0 & 3 & \-2 & \-2 & 2 & 1 & 0 & 0 & 1 & 1 & 2 & 0 & 0 & \-2 & \-1 & 1 \cr
 \-1 & 0 & 2 & \-1 & 1 & 2 & 0 & 2 & \-1 & \-1 & 1 & \-1 & 1 & \-1 & 2 & \-3 & 0 & \-3 & \-2 & 2 & \-2 & 2 & 3 & 0 & 2 & 3 & \-1 & 1 & \-2 & \-2 & \-1 & \-1 & 0 & 3 & \-1 & 3 & 0 & 0 & 0 & \-3 & 0 & \-2 & \-1 & 0 & 0 & 1 & 1 & 0 & \-1 & \-1 & 0 & 1 & 2 & \-1 & 1 \cr
 1 & 1 & \-1 & 1 & \-1 & 2 & \-1 & 0 & \-3 & 1 & 0 & \-1 & 0 & \-1 & \-2 & \-2 & 1 & 2 & 0 & 0 & 0 & 0 & 0 & \-1 & 1 & \-1 & 1 & 0 & \-2 & \-1 & 0 & \-2 & \-1 & 1 & 1 & 0 & 1 & 0 & \-1 & \-2 & 1 & 1 & 1 & 0 & \-2 & \-1 & 0 & \-3 & \-2 & 3 & \-1 & 0 & 1 & \-1 & 0 \cr
 0 & 0 & 1 & 2 & 1 & \-1 & 0 & 3 & \-2 & 1 & 0 & 0 & 2 & 1 & 2 & 0 & 0 & 1 & 1 & 1 & 1 & 0 & 1 & \-1 & 2 & 0 & 1 & \-1 & 3 & 0 & \-4 & 0 & \-2 & 1 & 0 & 0 & 1 & 1 & \-2 & 0 & 0 & 0 & \-2 & 1 & \-1 & 1 & 2 & 1 & \-3 & 0 & \-1 & 0 & \-1 & \-1 & 0 \cr
 \-1 & 0 & \-2 & \-1 & 0 & \-1 & \-2 & \-1 & 0 & \-1 & \-2 & 0 & 1 & \-1 & 2 & 0 & 1 & \-2 & 0 & \-1 & \-1 & 0 & 1 & 2 & 0 & 2 & \-4 & 1 & \-3 & \-1 & 0 & 1 & 1 & \-2 & 2 & \-2 & \-1 & 0 & 2 & 0 & \-1 & \-3 & \-3 & 1 & 2 & \-2 & \-2 & 0 & 2 & 2 & 2 & 1 & \-3 & 1 & \-1 \cr
 \-3 & \-1 & \-1 & 2 & 1 & \-3 & 1 & 2 & \-1 & \-1 & 1 & \-1 & 0 & 1 & 0 & 0 & \-1 & 2 & \-1 & \-1 & \-1 & 1 & 2 & \-1 & \-1 & \-1 & 2 & 1 & 1 & \-1 & \-2 & 1 & 1 & \-2 & 1 & 1 & 0 & 1 & 0 & \-1 & 1 & \-2 & \-2 & \-1 & 2 & 1 & \-1 & 2 & \-1 & \-1 & 0 & 0 & \-2 & \-2 & \-1 \cr
 \-3 & 0 & \-2 & 1 & 4 & 2 & \-1 & 0 & 1 & \-3 & 0 & \-2 & 0 & 1 & 2 & 0 & 0 & \-2 & 0 & 0 & \-1 & 1 & 3 & \-1 & 0 & 3 & \-1 & 2 & \-3 & 0 & 0 & 1 & 3 & 2 & \-1 & 2 & \-1 & 1 & 1 & \-1 & 2 & \-1 & \-1 & \-2 & 3 & 0 & 1 & 1 & 0 & 2 & 0 & 2 & 2 & \-2 & \-2 \cr
 0 & \-2 & 3 & \-2 & \-1 & 4 & \-2 & \-1 & 1 & \-2 & 3 & 2 & 1 & 0 & 2 & 1 & 3 & \-1 & \-1 & 1 & \-1 & 1 & 2 & 2 & 1 & 1 & \-1 & 1 & \-1 & \-3 & 2 & \-1 & 1 & 3 & \-3 & 1 & \-1 & \-1 & 2 & \-1 & 0 & 0 & 1 & \-1 & \-1 & 0 & 2 & \-1 & \-1 & \-1 & 1 & \-1 & 2 & \-1 & 2 \cr
 1 & \-1 & 2 & \-1 & \-1 & 2 & 0 & 0 & 1 & 1 & 0 & 2 & \-1 & 0 & 0 & 3 & 0 & 1 & \-1 & 0 & 2 & \-2 & 0 & \-1 & 0 & \-2 & 2 & 0 & 0 & 1 & 0 & \-1 & 0 & 2 & \-2 & 2 & 0 & 1 & 0 & 0 & 0 & 1 & 2 & 0 & 1 & 0 & 1 & 0 & 0 & \-1 & 0 & \-1 & 2 & 0 & 0 \cr
 0 & 0 & 2 & \-2 & \-1 & 0 & \-4 & 0 & \-2 & 2 & \-1 & \-1 & 2 & \-2 & 2 & 1 & 2 & \-3 & 1 & 1 & \-1 & \-1 & \-1 & 0 & 1 & 3 & \-2 & \-2 & 0 & \-2 & 0 & 1 & \-1 & 1 & 0 & 0 & 0 & \-2 & 1 & 0 & \-2 & 0 & \-1 & 3 & \-1 & \-2 & 1 & 1 & 1 & 0 & 0 & 1 & \-2 & 2 & 1 \cr
 0 & 1 & 1 & \-2 & \-1 & 1 & 1 & 1 & 2 & 0 & 1 & 1 & 0 & \-2 & 0 & \-1 & 0 & \-1 & \-1 & 0 & 1 & \-2 & 1 & 3 & \-1 & \-2 & 0 & 1 & \-1 & 1 & 2 & \-1 & 4 & \-3 & 0 & 2 & \-1 & 1 & 2 & \-3 & \-2 & \-1 & 2 & 0 & 0 & \-1 & 0 & \-3 & 3 & \-1 & 1 & 1 & 3 & \-1 & 0 \cr
 \-2 & 1 & \-3 & 2 & 0 & \-3 & 1 & 0 & \-3 & 2 & \-4 & \-3 & \-2 & \-1 & \-2 & \-1 & \-2 & 1 & 0 & \-2 & \-1 & 0 & 0 & \-2 & 1 & \-1 & 1 & 0 & 0 & 1 & \-2 & 1 & \-2 & \-3 & 5 & 0 & 0 & 2 & 0 & 1 & 0 & \-2 & \-1 & 1 & 3 & \-1 & \-3 & 2 & 2 & 2 & \-2 & 2 & \-2 & 1 & \-1 \cr
 0 & \-2 & \-2 & 1 & \-1 & \-2 & 1 & \-1 & 3 & \-2 & 1 & 1 & \-1 & \-1 & 1 & 1 & \-1 & 1 & \-1 & \-2 & 0 & 0 & 1 & 0 & \-3 & \-2 & \-1 & 3 & \-1 & 0 & 0 & \-1 & 1 & \-2 & 1 & \-2 & 0 & \-1 & 0 & 1 & 0 & \-2 & 0 & 0 & 2 & 0 & \-1 & 0 & 1 & \-1 & 2 & \-2 & \-1 & 1 & \-1 \cr
 \-1 & \-2 & \-1 & 4 & 0 & 4 & \-2 & \-1 & 0 & \-1 & \-1 & 0 & 1 & 0 & 0 & 2 & 0 & 0 & 1 & \-3 & 1 & 0 & 1 & \-2 & 2 & 0 & \-1 & 3 & 1 & \-2 & \-1 & 0 & 0 & 3 & \-1 & \-1 & \-1 & \-3 & 0 & 1 & 4 & 1 & 0 & 0 & 2 & \-1 & 1 & 1 & \-1 & 4 & \-1 & \-1 & 0 & \-2 & 0 \cr
 \-3 & \-2 & 1 & \-3 & 0 & 1 & \-2 & \-2 & 2 & \-3 & 2 & 0 & \-3 & 1 & 0 & 2 & 2 & \-1 & \-1 & \-1 & \-3 & 2 & 0 & 1 & 0 & 0 & 0 & 1 & \-1 & \-1 & 4 & 1 & 0 & 0 & \-1 & 2 & \-3 & 0 & 4 & 1 & 0 & \-1 & 1 & \-4 & 2 & 0 & \-1 & 3 & 0 & \-1 & 1 & \-1 & 1 & 1 & 2 \cr
 0 & 0 & 0 & 0 & 1 & 0 & 0 & 0 & 2 & 0 & 1 & \-1 & 1 & 1 & 0 & \-1 & 1 & 1 & \-1 & 1 & 0 & 0 & 0 & \-1 & \-1 & \-2 & 1 & \-1 & 2 & 2 & 0 & 1 & 1 & \-1 & \-1 & 1 & \-1 & \-1 & 2 & 1 & \-1 & 1 & 1 & 0 & \-1 & 0 & 1 & 0 & 1 & 1 & 1 & 0 & 1 & 0 & 0 \cr
 1 & 1 & \-1 & \-1 & 4 & \-3 & 2 & \-1 & 1 & 1 & \-1 & 1 & \-1 & 2 & 2 & 1 & \-3 & 0 & 1 & 0 & 1 & 2 & \-2 & \-1 & \-1 & 3 & \-2 & 0 & 2 & 2 & \-2 & 1 & \-1 & 2 & 2 & 1 & 1 & 2 & \-4 & 1 & \-2 & 0 & 0 & \-2 & 0 & 3 & 1 & 1 & \-2 & \-3 & 2 & 1 & 0 & 2 & 0 \cr
 2 & \-1 & 2 & 0 & \-3 & 1 & 0 & 1 & \-1 & \-2 & 2 & 0 & \-1 & 0 & \-1 & \-1 & 0 & 0 & 0 & 2 & \-1 & 0 & 0 & 1 & 1 & 0 & 1 & \-2 & \-2 & \-1 & 1 & \-2 & \-2 & 2 & \-4 & \-3 & 1 & 0 & \-1 & \-1 & 2 & 1 & \-1 & 2 & \-1 & 1 & 0 & 0 & \-3 & 0 & \-2 & \-2 & 0 & \-1 & \-1 \cr}
 $$}

\bigskip\noindent
The  prime $23$ is inert in the imaginary quadratic orders ${\cal O}$ of discriminant $-3,-4,-8$, $-12$, $-16$, $-27$ and $-163$.
For each order ${\cal O}$ there is a rational CM-point   on $X_{\rm ns}^+(23)$,  corresponding to an elliptic curve with complex multiplication by~${\cal O}$.  The CM-points have been computed numerically.  They are the only rational points $(x_1:x_2:x_3:x_4:x_5:x_6:x_7:x_8:x_9:x_{10}:x_{11}:x_{12}:x_{13})$ with  $x_i\in\ZZ$ satisfying $|x_i|\le 10\,000$.

\medskip
\hbox{\vbox{\hsize 2in\noindent {\bf Table 4.} CM-points on $X_{\rm ns}^+(23)$.\medskip\vbox{\offinterlineskip
\hrule\halign{&\vrule#&\strut\quad\hfil#\quad\cr
height3pt
&\omit&&\omit&\cr
&discriminant&&CM-point&\cr
height3pt
&\omit&&\omit&\cr
\noalign{\hrule}
height2pt
&\omit&&\omit&\cr
&$-3$&&$(-3:4:0:1:0:6:-1:6:-6:-6:0:-6:-12) $&\cr
&$-4$&&$(1:-2:0:-2:-1:0:1:-2:-1:0:-1:0:0) $&\cr
&$-8$&&$(3:13:-19:-4:16:8:-11:10:1:-7:-12:18:-5)$&\cr
&$-12$&&$(-15:4:-20:-3:12:6:9:-4:18:12:14:2:2) $&\cr
&$-16$&&$(3:-10:4:-4:-7:8:-11:10:1:16:11:18:18) $&\cr
&$-27$&&$ (0:1:0:1:0:0:-1:0:0:0:0:0:0) $&\cr
&$-163$&&$(0:-1:0:-1:0:-2:1:-2:-4:0:4:-2:2) $&\cr
height3pt
&\omit&&\omit&\cr
}\hrule}}}

\bibliography

\item{[\AL]}{Atkin, A. and  Li, W.: Twists of Newforms and Pseudo-Eigenvalues of $W$-Operators, {\sl Inventiones math.} {\bf 48} (1978), 221--243.
\item{[\BD]} Balakrishnan, J., Dogra, N.,  M\"uller, J.S., Tuitman, J. and  Vonk, J.:
Explicit Chabauty-Kim for the Split Cartan Modular Curve of Level $13$, arXiv:1711.05846, November 2017.
\item{[\BC]} Banwait, B. and Cremona, J.: Tetrahedral elliptic curves and the local-global principle for isogenies},
{\sl Algebra and Number Theory}, {\bf 8} (2014) 1201--1229.
\item{[\BA]}  Baran, B.: Normalizers of non-split Cartan subgroups, modular curves and the class number one problem, {\sl Journal of Number Theory} {\bf 130} (2010) 2753--2772.
\item{[\BB]}  Baran, B.: An exceptional isomorphism between modular curves of level 13, {\sl Journal of Number Theory} {\bf 145} (2014) 273--300.
\item{[\BPR]} Bilu, Y.,  Parent, P. and  Rebolledo, M.:   Rational points on $X_0^+ (p^r)$, {\sl Annales de l'institut Fourier}, {\bf 63} (2013), 957--984
\item{[\BK]} Birch, B. and Kuyk, W.: {\sl Modular Functions of One Variable IV},
Springer Lecture Notes {\bf 476}, Springer-Verlag 1972.
\item{[\BU]} Bump, D.: {\sl Automorphic Forms and Representations}, Cambridge Studies in Advanced Mathematics {\bf 55}, CUP 1998.
\item{[\KCON]} Conrad, K.: On Weil's proof of the bound for Kloosterman sums, {\sl J. Number Theory} {\bf 97} (2002), 439--446.
\item{[\DSED]} De Smit, B. and  Edixhoven, B.:  Sur un r\'esultat d'Imin Chen, {\sl Math. Res. Lett.} {\bf 7} (2000) 147--153.
\item{[\DGGS]} Dose, V., Fern\'andez, J., Gonz\'alez, J. and Schoof, R.: The automorphism group of the non-split Cartan modular curve of level $11$, {\sl Journal of Algebra} {\bf 417} (2014), 95--102.
\item{[\DOSE]} Dose, V.: On the automorphisms of the non-split Cartan modular curves of prime level, {\sl Nagoya Math. J.} {\bf 224} (2016),  74--92.
\item{[\DMS]} Dose, V., Mercuri, P. and  Stirpe, C.: Cartan modular curves of level $13$, to appear in {\sl Journal of Number Theory}.
\item{[\LAN]} Lang, S.: {\sl Algebra}, Graduate Texts in Math {\bf 73}, Springer-Verlag 2002.
\item{[\BAO]}  Le Hung, Bao V.: Modularity of some elliptic curves over totally real fields, preprint 2013, 
arXiv:1309.4134
\item{[\LI]}  Ligozat, G.: Courbes modulaires de niveau 11, pp. 149--237 in: {\sl Modular functions of one variable, V },
Lecture Notes in Math. {\bf 601} Springer-Verlag, Berlin, 1977. 
\item{[\MDB]} The modular forms data base LMFDB, {\tt http://www.lmfdb.org}.
\item{[\MAG]} Magma Computational Algebra System. {\tt http://magma.maths.usyd.edu.au/magma/}
\item{[\MAI]} Mazur, B.: Rational isogenies of prime degree, {\sl Invent. Math.} {\bf 44} (1978) 129--162.
\item{[\ME]} Mercuri, P.: {\sl Rational Points
on Modular Curves}, PhD thesis, Universit\`a di Roma ``La Sapienza" 2013/2014.
\item{[\MER]} Mercuri, P.: Equations and rational points of the modular curves $X_0^+(p)$, {\sl The Ramanujan Journal}, 
doi: 10.1007/s11139-017-9925-2, to appear.
\item{[\SAG]} Sagemath:  {\tt http://www.sagemath.org}.
\item{[\SD]} Saint-Donat, B.: On Petri's analysis of the linear system of quadrics through a canonical curve, {\sl Math. Ann.} {\bf 206} (1973), 157--175.
\item{[\SE]} Serre, J.-P.: Propri\'et\'es galoisiennes des points d'ordre fini des courbes elliptiques, {\sl Invent. Math.} {\bf 15} (1972)  259--331.
\item{[\ST]} Stein, W.: The modular forms database, {\tt http://wstein.org/Tables/tables.html}.
\item{[\WAL]} Waldspurger, J.-L.: Quelques propri\'et\'es arithm\'etiques de certaines formes automorphes sur ${\rm GL}(2)$,
{\sl Compositio Math.} {\bf 54} (1985), 121--171.

\bye